\documentclass{amsart}

\usepackage{amssymb,latexsym,amsmath,graphics,setspace,amscd,verbatim}
\usepackage{epsfig,epsf,amsthm,amsfonts}
\usepackage[active]{srcltx}

\theoremstyle{plain}
\newtheorem{theo}{Theorem}[section]
\newtheorem{prop}[theo]{Proposition}
\newtheorem{lemma}[theo]{Lemma}

\theoremstyle{definition}
\newtheorem{defi}[theo]{Definition}
\newtheorem{remark}[theo]{Remark}

\newcommand{\R}{\mathbb{R}}
\newcommand{\Z}{\mathbb{Z}}
\newcommand{\A}{\mathcal{A}}

\newcommand{\C}{\mathbb{C}}

\newcommand{\D}{\mathbb{D}}

\renewcommand{\P}{{\mathcal{P}}}

\newcommand{\U}{{\mathcal{U}}}

\newcommand{\interior}[1]{\mathring{{#1}}}

\newcommand{\supp}{\text{supp}}

\renewcommand{\sl}{\text{sl}}

\newcommand{\dist}{\text{dist}}
\newcommand{\link}{\text{link}}

\begin{document}

\title
[Reeb flows and applications to Finsler geodesic flows]{Global properties of tight Reeb flows with applications to Finsler geodesic flows on $S^2$}
\author{Umberto L. Hryniewicz}
\address[Umberto L. Hryniewicz]{Institute for Advanced Study/Departamento de Matem\'atica Aplicada -- Universidade Federal do Rio de Janeiro}
\email{umbertolh@math.ias.edu/umberto@labma.ufrj.br}
\author{Pedro A. S. Salom\~ao}
\address[Pedro A. S. Salom\~ao]{Departamento de Matem\'atica, Instituto de Matem\'atica e Estat\'istica -- Universidade de S\~ao Paulo}
\email{psalomao@ime.usp.br}

\begin{abstract}
We show that if a Finsler metric on $S^2$ with reversibility $r$ has flag curvatures $K$ satisfying $(\frac{r}{r+1})^2 <K \leq 1$, then closed geodesics with specific contact-topological properties cannot exist, in particular there are no closed geodesics with precisely one transverse self-intersection point. This is a special case of a more general phenomenon, and other closed geodesics with many self-intersections are also excluded. We provide examples of Randers type, obtained by suitably modifying the metrics constructed by Katok~\cite{katok}, proving that this pinching condition is sharp. Our methods are borrowed from the theory of pseudo-holomorphic curves in symplectizations. Finally, we study global dynamical aspects of $3$-dimensional energy levels $C^2$-close to $S^3$.
\end{abstract}

\maketitle

\section{Introduction and Main Results}

Classical and recent results show that pinching conditions on the curvatures of a Riemannian metric force the geodesic flow to present specific global behavior, usually encoded in geometric-topological and dynamical properties of closed geodesics.  The interest in such phenomena can be traced back to Poincar\'e~\cite{po} and Birkhoff~\cite{bi} where, among many other topics, the geodesic flow on positively curved surfaces was studied.

In the 1970s and 1980s this subject again received much attention. For example, in the work of Thorbergsson~\cite{thor}, Ballmann, Thorbergsson and Ziller \cite{btz1,btz2,btz3} and Klingenberg~\cite{kli} one finds many results relating pinching conditions on the curvatures to the existence (or non-existence) of closed geodesics with various topological and dynamical properties. Let us recall two theorems along these lines proved around the same time. As usual, a Riemannian metric is called $\delta$-pinched if all sectional curvatures $K$ satisfy $\delta \leq K\leq 1$.

\begin{theo}[Ballmann~\cite{ball}]\label{ball_theo}
Given $k\geq 1$ and $\epsilon>0$, there exists $\delta<1$ such that every prime closed geodesic of a $\delta$-pinched metric on $S^2$ is either a simple curve of length in $[2\pi-\epsilon,2\pi+\epsilon]$, or has at least $k$ self-intersections and length $>\epsilon^{-1}$.
\end{theo}

\begin{theo}[Bangert~\cite{bangert}]\label{bang_theo}
For every $\epsilon>0$ there exists $\delta<1$ such that the length $l$ of every prime closed geodesic of a $\delta$-pinched metric on $S^n$ satisfies either $l \in [2\pi-\epsilon,2\pi+\epsilon]$ or $l>\epsilon^{-1}$.
\end{theo}

Both theorems are of a perturbative nature and exhibit a ``short-long'' dichotomy for prime closed geodesics: if the metric is sufficiently pinched then their lengths are either close to the lengths in the round case (short), or arbitrarily large. This is surprising since one could try to imagine a sequence of metrics converging in $C^2$ to the round sphere admitting prime closed geodesics with lengths close to $2k\pi$, for some $k\geq 2$. However, this does not happen.

Motivated by the above statements one might consider the following questions in the more general framework of Finsler metrics on $S^2$, or even in broader classes of Hamiltonian systems.
\begin{enumerate}
  \item[a)] How much can we relax the pinching of the flag curvatures of a (possibly non-reversible) Finsler metric on the $2$-sphere and still keep some kind of dichotomy similar to that in Theorem~\ref{ball_theo}?
  \item[b)] Can the ``short-long'' length dichotomy in Theorem~\ref{bang_theo} be generalized to a ``low-high'' action dichotomy on a broader class of Hamiltonian systems? If so, is there any additional topological information that can be extracted in low dimensions, in a way similar to Theorem~\ref{ball_theo}?
\end{enumerate}

Ballmann, Thorbergsson and Ziller~\cite{btz3} observe that a $\delta$-pinched Riemannian metric on the $2$-sphere, with $\delta>1/4$, does not admit a closed geodesic with precisely one self-intersection. The proof is an immediate application of two well-known comparison theorems.  To be more precise,  the pinching condition $0 < \delta<K\leq 1$ implies Klingenberg's estimate for the injectivity radius $\text{inj}(p) \geq \pi,\forall p\in S^2$. Since a closed geodesic $\gamma$ with exactly one self-intersection point is the union of two loops, its lentgh $l$ must satisfy $l\geq 4\pi$. On the other hand, since $\gamma$ is also a convex geodesic polygon, we have from Toponogov's theorem the estimate $l \leq 2\pi / \sqrt{\delta}< 4\pi$ if $\delta>1/4$. These two inequalities on the length $l$ imply that such closed geodesic cannot exist. In this case, closed geodesics are either simple with length $\leq2\pi/\sqrt{\delta}$, or have at least two self-intersections and length $\geq 6\pi$. This may be thought of as a simplest answer to a) in the Riemannian case, but maybe other pinching conditions will rule out other types of geodesics.

We use the theory of pseudo-holomorphic curves in symplectizations developed by H. Hofer, K. Wysocki and E. Zehnder as an alternative to comparison theorems. These methods reveal a more general phenomenon, in fact, under a certain pinching condition on the flag curvatures (cf. Theorem~\ref{main1}) there exists a larger class of immersed curves that cannot be realized as closed geodesics. This class includes curves with precisely one transverse self-intersection, but also many other curves with an arbitrarily large number of self-intersections. Then we exhibit examples of Randers type showing that the above mentioned pinching condition is optimal.

Finally, we quickly address b). It follows trivially from the method of Bangert~\cite{bangert} that a low-high action dichotomy holds for convex energy levels in $\R^{2n}$ which are $C^2$-close to $S^{2n-1}$. In the case $n=2$ we study the linking number between high- and low-action orbits.

\subsection{Main results}


We consider a weakened version of the notion of flat knot types discussed in~\cite{angenent}, which relates to V. I. Arnold's ${\bf{\rm J^+}}$-theory of plane curves described in~\cite{arnold}.

\begin{defi}\label{weak_flat_knot}
Consider the set $\mathcal B$ of $C^2$-immersions $\gamma : S^1 \to S^2$ such that all self-intersections are either transverse or negative tangencies, i.e., if $\gamma(t_0) = \gamma(t_1)$ and $t_0 \not= t_1$ then $\dot\gamma(t_1) \not\in \R^+\dot\gamma(t_0)$. We say two curves $\gamma_0,\gamma_1 \in \mathcal B$ are equivalent if they are homotopic through curves in $\mathcal B$. A weak flat knot type is an equivalence class of curves in $\mathcal B$.
\end{defi}


If we fix a Finsler metric $F$ on $S^2$ then the unit sphere bundle $SS^2 := \{v \in TS^2 \mid F(v)=1 \}$ admits a contact form $\alpha_F$ given by the pull-back of the tautological $1$-form of $T^*S^2$ via the associated Legendre transform. Any weak flat knot type of some $\gamma \in \mathcal B$ singles out a transverse knot type in the contact manifold $(SS^2,\ker\alpha_F)$ determined by the knot $\dot\gamma/F(\dot\gamma)$. In particular, topological and contact invariants of transverse knots, like the self-linking number (cf. \S~\ref{sl_defn} below), induce invariants of weak flat knot types on $S^2$. Note that any prime closed geodesic of a Finsler metric on $S^2$ represents a weak flat knot type\footnote{In the reversible case negative self-tangencies of closed prime geodesics never happen, so one gets a so-called flat knot type as defined in~\cite{angenent}. However, negative self-tangencies could appear in the non-reversible case.}.

To give a computable concrete example, consider the weak flat knot type $k_8$ of an ``eight-like curve'' having precisely one self-intersection point which is transverse. The proof of the following lemma is found in \S~\ref{par_topol}.

\begin{lemma}\label{lema_topol}
Let the $C^2$-immersion $c : S^1 \to S^2$ represent the weak flat knot type $k_8$, and let $F$ be any Finsler metric on $S^2$. Then the curve $\dot c/F(\dot c)$ in the unit sphere bundle is unknotted and has self-linking number $-1$.
\end{lemma}


Before stating our main result we need to recall the notion of reversibility of a Finsler metric $F$, defined by Rademacher~\cite{rad1} as
\begin{equation}\label{rev_def}
  r := \max \{ F(-v) \mid F(v)=1 \} \geq 1.
\end{equation}
It equals $1$ exactly when $F(v) = F(-v) \ \forall v$, and $F$ is called reversible in this case. The notion of reversibility is an essential ingredient in Rademacher's proof of his sphere theorem for Finsler metrics. 

\begin{figure}
  \includegraphics[width=300\unitlength]{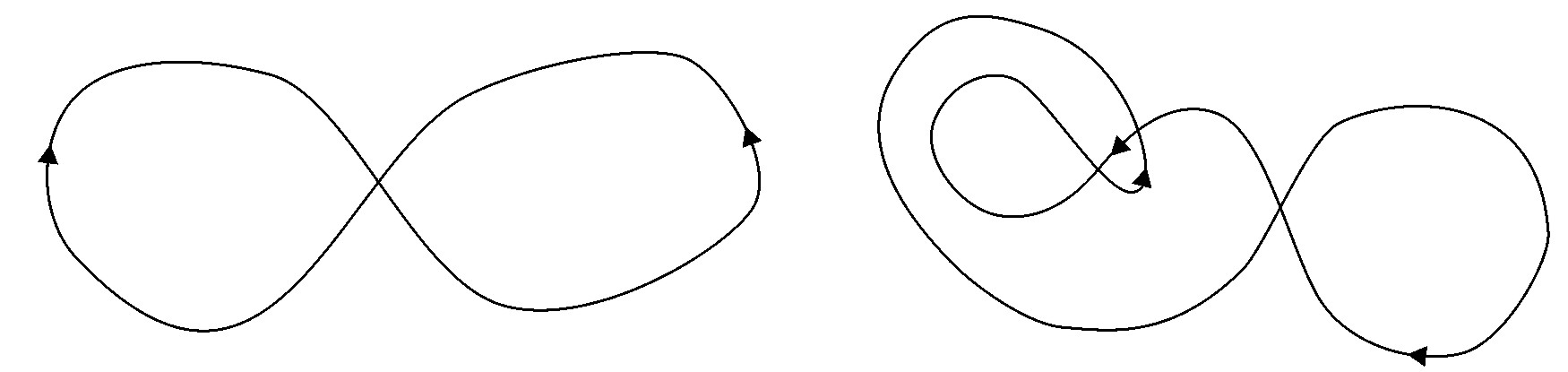}
  \caption{The weak flat knot type $k_8$.}
\end{figure}

\begin{theo}\label{main1}
The following assertions hold.
\begin{itemize}
  \item[(i)] Let $F$ be a Finsler metric on $S^2$ with reversibility $r$. If all flag curvatures $K$ satisfy
      \begin{equation}\label{hypothesis1}
        \left( \frac{r}{1+r} \right)^2 < K \leq 1
      \end{equation}
      then no prime closed geodesic $\gamma$ represents the weak flat knot type $k_8$.
  \item[(ii)] Statement (i) is optimal in the following sense: for every choice of $r\geq 1$ and~$0<\delta<(r/(r+1))^2$ there exists a Finsler metric on $S^2$ with reversibility $r$ and $\delta$-pinched flag curvatures admitting closed geodesics with precisely one transverse self-intersection.
\end{itemize}
\end{theo}

We stress the fact that the proof of part (i) in Theorem~\ref{main1} does not make use of any version of Toponogov's theorem for Finsler geometry. In fact, our method seems to be an alternative tool in those cases where such a comparison theorem may not be effective.

Note that there exist immersions $\gamma \in \mathcal B$ representing $k_8$ with an arbitrarily large number of self-intersections, see Figure 1 for an example with 3 self-intersections. Such immersions cannot be realized as a closed geodesic under the pinching condition~\eqref{hypothesis1}.

To prove assertion (ii) we modify the metrics of Katok~\cite{katok}. The examples are Randers metrics given by suitably chosen Zermelo navigation data on surfaces of revolution in $\R^3$, see \S~\ref{examples} for the detailed construction. Assertion (i) is proved by an application of pseudo-holomorphic curve theory in symplectizations, as introduced by Hofer in~\cite{93}, developed by Hofer, Wysocki and Zehnder during the 1990s, and later by many other authors. The arguments are based on a dynamical characterization of the tight $3$-sphere from~\cite{hryn1}, extending earlier results from~\cite{char1,char2}. For an outline of the proof we refer to \S~\ref{sketch_main1} below.

Our second result relates to question b). Denote by $\textbf{Conv}(2m)$ the set of closed and strictly convex hypersurfaces of class $C^2$ in $\R^{2m}$, equipped with the $C^2$-topology. $\R^{2m}$ is endowed with its standard symplectic structure $\omega_0$, and each $S \in \textbf{Conv}(2m)$ is oriented as the boundary of the bounded connected component of $\R^{2m}\setminus S$. A closed characteristic on $S$ is a closed leaf of $TS^\bot \subset TS$, where $\bot$ denotes the $\omega_0$-symplectic orthogonal. These are precisely geometric images of closed Hamiltonian orbits, for any Hamiltonian realizing $S$ as a regular energy level. We denote by $\P(S)$ the set of closed characteristics and think its elements as first iterates of periodic orbits of a Hamiltonian system. The action of a given $P \in \P(S)$ is $\mathcal A(P) = \int_P \lambda_0$, where $\lambda_0 = \frac{1}{2} \sum qdp-pdq$ is the standard Liouville form.

These Hamiltonian systems generalize geodesic flows on $S^2$. In fact, as is well known, the geodesic flow of any Finsler metric on $S^2$ lifts to a Hamiltonian flow on a suitable star-shaped hypersurface in $\R^4$ via a double cover. This lifting procedure is nicely described in~\cite{HP}. In general, however, $S$ is not convex. If the metric is $C^2$-close to the round metric then $S$ belongs to $\textbf{Conv}(4)$ and is close to $S^3$. Action of an orbit on $S$ coincides with length of its projection on $S^2$. With this picture in mind we make the following statement.

\begin{theo}\label{main2}
Given $\epsilon>0$ there exists a neighborhood $\U_\epsilon$ of $S^{2m-1}$ in $\textbf{Conv}(2m)$ such that if $S \in \U_\epsilon$ then every $P\in \P(S)$ satisfies $\A(P) \in [\pi-\epsilon,\pi+\epsilon]$ (short orbits) or $\A(P) > \epsilon^{-1}$ (long orbits). In the case $m=2$, given any $k\geq 1$ there exists a neighborhood $\U_{\epsilon,k} \subset \U_\epsilon$ of $S^3$ such that if $S \in \U_{\epsilon,k}$ then $\link(P,P') \geq k$ whenever $P$ is short and $P'$ is long.
\end{theo}

The assertion about the high-low action dichotomy is a direct application of Proposition~\ref{propBangert} due to Bangert~\cite{bangert}, and in the case $m=2$ it is crucial to estimating the linking numbers. Obviously the short orbits are unknotted, have self-linking number $-1$ and their Conley-Zehnder indices belong to $\{3,4,5\}$. An analogous statement is true for energy levels $C^2$-close to irrational ellipsoids, except that the high-low action dichotomy is trivial in this case. For an idea of the proof see \S~\ref{sketch_main2} below.

\subsection{Outline of the main arguments}\label{outline}

For convenience of the reader we sketch some of the main steps in the proofs of our results.

\subsubsection{Non-existence of geodesics}\label{sketch_main1}

Here we briefly explain why assertion (i) of Theorem~\ref{main1} holds. For more details see \S~\ref{non-exist}.

A contact form $\lambda$ on a $3$-manifold is called non-degenerate if the spectrum of the linearized Poincar\'e map associated to any (prime) closed Reeb orbit does not contain roots of unity when restricted to the contact structure. According to~\cite{convex}, $\lambda$ is said to be dynamically convex if $c_1(\ker\lambda)$ vanishes and the Conley-Zehnder index of every contractible closed orbit of the associated Reeb flow is at least $3$. See \S~\ref{CZ-3d} for a definition of the index in $3$-dimensions.

In~\cite{char1,char2} it is proved that a closed connected tight contact $3$-manifold $M$ is the tight $3$-sphere if, and only if, the contact structure can be realized as the kernel of a dynamically convex non-degenerate contact form admitting an unknotted closed Reeb orbit $P$ with self-linking number $-1$ and Conley-Zehnder index $3$. In fact, they show that the given orbit bounds a disk-like global surface of section for the Reeb flow, but much more can be said: there is an open book decomposition of $M$ with disk-like pages and binding $P$, such that every page is a global surface of section. In particular, $M$ is homeomorphic to $S^3$. The following result from~\cite{hryn1} states that the restriction on the Conley-Zehnder index can be dropped.

\begin{theo}\label{secao_global_1}
Let $\lambda$ be a non-degenerate dynamically convex tight contact form on a closed connected $3$-manifold $M$. A closed Reeb orbit $P$ is the binding of an open book decomposition with disk-like pages which are global surfaces of section for the Reeb flow if, and only if, it is unknotted and has self-linking number $-1$. In particular, $M$ is homeomorphic to the 3-sphere when an orbit $P$ with these properties exists.
\end{theo}

The geodesic flow restricted to the unit sphere bundle of a Finsler metric $F$ coincides with the Reeb flow of the contact form $\alpha_F$, as explained before. Suppose $F$ is such a metric on $S^2$ satisfying~\eqref{hypothesis1}, and assume $\gamma$ is a closed geodesic whose lift $\dot\gamma$ is unknotted and has self-linking number $-1$ in $SS^2$. Theorem~\ref{secao_global_1} together with the statement below due to Harris and Paternain~\cite{HP} provides, in the bumpy case, a contradiction to the existence of $\gamma$ since the unit sphere bundle is not homeomorphic to the $3$-sphere. Thus the weak flat knot type $k_8$ can not be realized by a closed geodesic in view of Lemma~\ref{lema_topol}. The general case is discussed in Section~\ref{sec_main1}.

\begin{theo}[Harris and Paternain]\label{teo_hp}
If a Finsler metric $F$ with reversibility $r$ on $S^2$ is $\delta$-pinched, for some $\delta>(r/(r+1))^2$, then $\alpha_F$ is dynamically convex.
\end{theo}

The proof of Theorem~\ref{teo_hp} relies on Rademacher's estimate for the length of closed geodesic loops, see \S~\ref{proof_hp} below. As it will be clear,  the proof of Theorem~\ref{main1} shows that Harris-Paternain's pinching condition is sharp in the following sense: given $r\geq 1$ and  $0<\delta<(r/(r+1))^2$ there exists a $\delta$-pinched Finsler metric on $S^2$ with reversibility $r$, such that $\alpha_F$ is not dynamically convex.

\begin{remark}
As the reader may already have noticed, our argument will actually show a possibly stronger statement than that of Theorem~\ref{main1}. In fact, if we assume~\eqref{hypothesis1} then the transverse knot $\dot\gamma/F(\dot\gamma)$ associated to a closed prime geodesic $\gamma$ can not be unknotted and have self-linking number $-1$.
\end{remark}

\subsubsection{Convex energy levels $C^2$-close to $S^{2n-1}$}\label{sketch_main2}

The low-high action dichotomy in Theorem~\ref{main2} above is, of course, an immediate consequence of the non-trivial analysis from~\cite{bangert}.

Consider an unperturbed flow and a periodic orbit $P$ with prime period $T$ for which the linearized transverse Poincar\'e map is the identity. Then, roughly speaking, a prime closed orbit near $P$ of a perturbed flow either has period $\sim T$, or has a very large period. This is a particular instance of Proposition~\ref{propBangert} below which was extracted from~\cite{bangert}. The low-high action dichotomy is obtained when we take as the unperturbed flow the Reeb flow on $S^{2n-1}$ induced by the contact form $\lambda_0 = \frac{1}{2} \sum qdp-pdq$, since all orbits are periodic with the same period, and the transverse linearized Poincar\'e map is always the identity.

In the case $n=2$ we study the relation between orbits with low and high action. Analyzing specific global behavior of the ``round'' flow on $S^3$ we are able to conclude that a short orbit $P_s$ of the perturbed Reeb flow bounds a disk transverse to the flow. A long orbit $P_l$ either stays far from $P_s$, and thus links many times with it, or gets close to $P_s$ and again links many times since the linearized flow along $P_s$ rotates almost uniformly.

\section{Preliminaries}\label{prel_section}

This section is devoted to reviewing the definitions and facts necessary for the proofs that follow.

\subsection{The Conley-Zehnder index in $3$ dimensions}\label{CZ-3d}

The Conley-Zehnder index is an invariant of the linearized dynamics along closed Reeb orbits, which we now describe in the $3$-dimensional case.

Whenever $I \subset \R$ is a closed interval of length strictly less than $1/2$ satisfying $\partial I \cap \Z = \emptyset$, consider the integer $\hat\mu(I)$ defined by $\hat \mu(I) = 2k$ if $k \in I$, or $\hat\mu(I) = 2k+1$ if $I \subset (k,k+1)$. It can be extended to the set of all closed intervals of length strictly less than $1/2$ by $\hat\mu(I) = \lim_{\epsilon\to0^+} \hat\mu(I-\epsilon)$.

Let $\alpha$ be a contact form on the $3$-manifold $N$, inducing the contact structure $\xi = \ker \alpha$. Then $\xi$ becomes a symplectic vector bundle with the bilinear form $d\alpha|_\xi$. Suppose $x:\R\to N$ is a contractible periodic trajectory of the Reeb vector $R$ (uniquely defined by $i_R\alpha=1$ and $i_Rd\alpha=0$) of period $T>0$, and let $f:\D\to N$ be a map satisfying $f(e^{i2\pi t}) = x_T(t) := x(Tt)$. We can find a symplectic trivialization $f^*\xi \simeq \D\times \R^2$, which restricts to a trivialization $\Psi:x_T^*\xi \to S^1\times \R^2$. If $\phi_t$ is the Reeb flow then $d\phi_t$ preserves $\xi$, and we get a smooth path of symplectic $2\times 2$ matrices $t\in \R \mapsto \varphi(t) = \Psi_t \cdot d\phi_{Tt} \cdot \Psi_0^{-1}$, where $\Psi_t$ is restriction of $\Psi$ to $x_T^*\xi|_{t}$.

Given $w\in \R^2$, $w\not=0$, define $\Delta(w) := \frac{1}{2\pi}(\vartheta(1)-\vartheta(0))$, where $\vartheta(t)$ is a continuous lift of the argument of $\varphi(t)w$. Then consider the closed real interval $I := \{ \Delta(w) : w\in \R^2, \ w\not=0 \}$. It is easy to check that $I$ has length strictly less than $1/2$. 
Following~\cite{fols},  one can define
\begin{equation}\label{CZ_def}
  \mu_{CZ}(x,T) = \hat\mu(I).
\end{equation}
Note that $\mu_{CZ}(x,T) \geq 3$ if, and only if, $I \subset \{y\in\R : y > 1 \}$. Once $f$ is fixed, the above integer does not depend on the choice of symplectic trivialization $f^*\xi \simeq \D\times \R^2$. Nevertheless, the notation should indicate the dependence on the disk-map $f$, but $\mu_{CZ}(x,T)$ does not depend on $f$ when $c_1(\xi)$ vanishes, see~\cite{fols} for more details.

\subsection{The self-linking number}\label{sl_defn}

Let $(M,\xi)$ be a contact $3$-manifold, $L \subset M$ be a knot transverse to $\xi$, and let $\Sigma \hookrightarrow M$ be a Seifert surface for $L$, that is, $S$ is an orientable embedded connected compact surface $\Sigma \hookrightarrow M$ such that $L = \partial \Sigma$. Assume $\xi = \ker\lambda$ for some contact form $\lambda$. Since the bundle $\xi|_\Sigma$ carries the symplectic bilinear form $d\lambda$, there exists a smooth non-vanishing section $Z$ of $\xi|_\Sigma$ which can be used to slightly perturb $L$ to another transverse knot $L_\epsilon = \{ \exp_x (\epsilon Z_x) : x\in L \}$. Here $\exp$ is any exponential map. A choice of orientation for $\Sigma$ induces orientations of $L$ and of $L_\epsilon$. The {\it self-linking number} is defined as the oriented intersection number
\begin{equation}\label{defselflink0}
 sl(L,\Sigma) := L_\epsilon \cdot \Sigma \in \Z,
\end{equation}
where $M$ is oriented by $\lambda\wedge d\lambda$. It is independent of $\Sigma$ when $c_1(\xi) \in H^2(M)$ vanishes.

\subsection{Basics in Finsler geometry}\label{geomsetup}

We recall the basic definitions in Finsler geometry following~\cite{gri}. See also~\cite{rad,bao,vit}. The knowledgeable reader is encouraged to skip to \S~\ref{proof_hp}, and refer back only for the notation established here.

\subsubsection{Connections and curvatures}

Let $\pi:TM \to M$ be the tangent bundle of a manifold $M$, and denote $TM_0 := TM \setminus \{\text{zero section}\}$. Let $VTM$ be the vertical subbundle $\ker d\pi \subset TTM$, with fiber $V_vTM$ over $v\in TM$. $VTM_0$ denotes its restriction to $TM_0$. Whenever $(x^1,\dots,x^n)$ are coordinates on $M$ we have natural coordinates $$ (x^1,\dots,x^n,y^1,\dots,y^n) \simeq \sum_i y^i\partial_{x^i} $$ on $TM$. Thus, $\{\partial_{y^1},\dots,\partial_{y^n}\}$ is a local frame on $VTM$. On $TM$ we have a vector field defined in natural coordinates by $C = \sum_{i=1}^n y^i \partial_{y^i}$, and the almost tangent structure $\mathcal J$, which is the $VTM$-valued $1$-form on $TM$ defined locally by $\mathcal J = \sum_{i=1}^n dx^i \otimes \partial_{y^i}$. There is a canonical linear isomorphism $i_v: T_{\pi(v)}M \simeq V_vTM$, for any given $v \in TM$, defined by $i_v(w) = \left. \frac{d}{dt} \right|_{t=0} v+tw$. In natural coordinates: $i_v(w) = \sum_i w^i\partial_{y^i}$ if $w = \sum_i w^i\partial_{x^i}$. Thus $C_v = i_v(v)$.

A $TTM$-valued $1$-form $\Gamma$ on $TM_0$ satisfying
\begin{equation}\label{}
  \begin{array}{ccc}
    \Gamma^2 = I & \text{and} & \ker (\Gamma + I) = VTM_0
  \end{array}
\end{equation}
is a Grifone connection on $M$. In natural coordinates the equations $\Gamma \cdot \partial_{x^i} = \partial_{x^i} -2\Gamma^j_i \partial_{y^j}$, $\Gamma \cdot \partial_{y^i} = -\partial_{y^i}$ define the connection coefficients $\Gamma^j_i$ (we use Einstein summation convention). Considering the associated horizontal subbundle $HTM := \ker (\Gamma-I)$ we have a splitting
\begin{equation}\label{}
  TTM = VTM \oplus HTM
\end{equation}
and induced projections $P_V:TTM \to VTM$, $P_H:TTM \to HTM$. The isomorphisms $i_v^{-1}: V_vTM \simeq T_{\pi(v)}M$ and $d\pi: H_vTM \simeq T_{\pi(v)}M$ provide an isomorphism
\begin{equation}\label{iso_ttm}
  (i_v^{-1} \circ P_V,d\pi \circ P_H) : T_vTM \simeq T_{\pi(v)}M \oplus T_{\pi(v)}M
\end{equation}
when $v\not =0$. If $\zeta = \delta x^i\partial_{x^i} + \delta y^i\partial_{y^i} \in T_vTM$ then $$ \zeta \simeq \left( ( \delta y^i + \Gamma^i_k\delta x^k )\partial_{x^i} , \delta x^i \partial_{x^i} \right) \in T_{\pi(v)}M \oplus T_{\pi(v)}M $$ by the map~\eqref{iso_ttm}.

The curvature form of $\Gamma$ is the $VTM$-valued $2$-form on $TM_0$ defined by
\begin{equation}\label{}
  R(X,Y) = P_V([P_H(X),P_H(Y)])
\end{equation}
where $X,Y$ are vector fields on $TM_0$.

Later we will need to consider lifts of a Grifone connection $\Gamma$. These are linear connections $\nabla$ on $VTM$ satisfying
\begin{equation}\label{}
  HTM = \ker (X \mapsto \nabla_XC).
\end{equation}
$\nabla$ is said to be symmetric if $\nabla_X\mathcal J(Y) - \nabla_Y\mathcal J(X) = \mathcal J([X,Y])$ for arbitrary vector fields $X,Y$ on $TM_0$. If $\nabla$ has coefficients $\nabla_{\partial_{x^i}}\partial_{y^j} = \Gamma^k_{ij}\partial_{y^k}$ and $\nabla_{\partial_{y^i}}\partial_{y^j} = D^k_{ij}\partial_{y^k}$, in local natural coordinates, then this symmetry condition implies $D^k_{ij}=0$, $\Gamma^k_{ij}=\Gamma^k_{ji}$ and $y^j\Gamma^k_{ij} = \Gamma^k_i$.

The curvature tensor of $\nabla$ is
\begin{equation}\label{}
  \tilde R(X,Y)Z = \nabla_X\nabla_Y Z - \nabla_Y\nabla_X Z - \nabla_{[X,Y]}Z
\end{equation}
where $X,Y$ are vector fields on $TM_0$ and $Z$ is a section of $VTM_0$.
The curvature endomorphism of $\nabla$ in the direction of $v \in TM_0$ is the linear map $R^v : T_{\pi(v)}M \to T_{\pi(v)}M$ defined by
\begin{equation}\label{}
  R^v(u) = i_v^{-1}(\tilde R(v_h,u_h)i_v(v))
\end{equation}
where $u_h,v_h \in T_vTM$ are the (unique) horizontal lifts of $u,v$, respectively.

\subsubsection{Sprays and their geodesics}

Recall that a spray is a continuous vector field $S$ on $TM$, smooth on $TM_0$, satisfying equations
\begin{equation}\label{spray}
  \begin{array}{cc}
    i_S\mathcal J=C, & L_CS=S.
  \end{array}
\end{equation}
In local natural coordinates one can write $S = y^i\partial_{x^i} - 2G^i(x,y)\partial_{y^i}$. The $G^i$ will be referred to as the spray coefficients, and they satisfy $G^i(x,ty) = t^2G^i(x,y)$.

Every spray $S$ defines a Grifone connection by $\Gamma_S := -L_S\mathcal J$. It follows from~\eqref{spray} that $L_C\Gamma_S=0$ and that $\Gamma_S$ is symmetric: if $\Gamma^j_i$ are the connection coefficients in natural coordinates then $\partial_{y^k}\Gamma^j_i = \partial_{y^i} \Gamma^j_k$ for every $i,j,k$. Moreover, $\Gamma^j_i = \partial_{y^i}G^j$.

The set of symmetric lifts of $\Gamma_S$ is non-empty: the Berwald connection of $S$ is given in natural coordinates by
\[
  \begin{array}{cccc}
    \nabla_{\partial_{y^i}}\partial_{y^j} = 0, & \nabla_{\partial_{x^i}}\partial_{y^j} = \Gamma^k_{ij}\partial_{y^k} & \text{where} & \Gamma^k_{ij} = \partial_{y^j}\Gamma^k_i = \frac{\partial^2 G^k}{\partial y^j \partial y^i}.
  \end{array}
\]

If $\zeta(t) \subset TM_0$ is an integral curve of $S$ then $\zeta = \dot\gamma$ where $\gamma = \pi \circ \zeta$. This follows from~\eqref{spray}. A curve $\gamma(t)$ on $M$ is called a geodesic if $\dot\gamma(t)$ is an integral curve of $S$, where the condition $\dot\gamma(0) \not=0$ is implicit.

Using $\Gamma_S$ one defines the covariant derivative of a vector field $V(t)$ along a geodesic $\gamma(t)$. Namely, $\tilde V(t) := i_{\dot\gamma(t)}(V(t))$ defines a (vertical) vector field along the integral curve $\dot\gamma$ of $S$ and, hence, the Lie derivative $L_S\tilde V$ is well-defined. We set
\begin{equation}\label{cov_der}
  \frac{D_{\gamma}V}{dt} := i_{\dot\gamma(t)}^{-1}(P_V(L_S\tilde V)).
\end{equation}
In natural coordinates, if $V = V^i\partial_{x^i}$ and $\gamma^i(t) = x^i\circ\gamma$ then
\begin{equation}\label{cov_der_coord}
  \frac{D_{\gamma}V}{dt} = \left( \dot V^i + \Gamma^i_kV^k \right) \partial_{x^i}
\end{equation}
where the $\Gamma^i_k$ are evaluated at $(\gamma^1,\dots,\gamma^n,\dot\gamma^1,\dots,\dot\gamma^n)$. Note the drastic difference to the Riemannian case, where the $\Gamma^i_k$ do not depend on $\dot\gamma^1,\dots,\dot\gamma^n$. Thus, in the more general present situation, covariant differentiation along arbitrary curves may not be defined only in terms of the spray, and other choices must be made.

Parallel transport $P_t : T_{\gamma(0)}M \to T_{\gamma(t)}M$ along a geodesic $\gamma$ is defined by $P_t ( V_0) := V(t)$, where $V_0 \in T_{\gamma(0)}M$, $\frac{D_{\gamma}V}{dt}=0$ and $V(0)=V_0$.

\begin{lemma}\label{useful_lemma}
Let $\nabla$ be a symmetric lift of $\Gamma_S$. Then $$ -R^v(u) = i_v^{-1}(R(S,u_h)) \text{ for any } u\in T_{\pi(v)}M $$ where $u_h \in H_vTM$ is the horizontal lift of $u$ ($d\pi \cdot u_h = u$) and $R$ is the curvature form of $\Gamma_S$. In particular, $R^v$ is independent of the choice of the symmetric lift.
\end{lemma}

For completeness, and convenience of the reader, we include a proof of the above well-known standard fact in the appendix.

\subsubsection{The case of Finsler manifolds}\label{casefinsler}

A Finsler metric on a mani\-fold $M$ is a continuous function $F:TM \to [0,+\infty)$, smooth on $TM_0:=TM\setminus \{\text{zero section}\}$ satisfying:
\begin{itemize}
\item[(i)] $F(tv)=tF(v),\forall  v\in TM,\forall t>0$. $F$ is said to be positively homogeneous of degree $1$.
\item[(ii)] For each $v\in TM_0$ the quadratic form $g_v: T_{\pi(v)}M \times T_{\pi(v)}M \to \R$ given by
\begin{equation}\label{quad_form}
  g_v(w_1,w_2) = \frac{1}{2}\left. \frac{\partial^2}{\partial s\partial t} \right|_{s=t=0} F^2(v+sw_1 + tw_2).
\end{equation}
is positive definite. This is called the convexity condition.
\end{itemize}

The Legendre transform $\mathcal L : T^*M \to TM$ is the fiber-preserving homeomorphism defined in the following manner: given $\lambda \in T^*_pM$ then
\begin{equation}\label{legendre}
  \begin{array}{ccc}
    \mathcal L(\lambda) := v^* \in T_pM & \text{where} & \lambda \cdot v^* - \frac{1}{2} F^2(v^*) = \sup_{v\in T_pM} \lambda \cdot v - \frac{1}{2} F^2(v).
  \end{array}
\end{equation}
Thus $\mathcal L(t\lambda) = t\mathcal L(\lambda)$, $\forall t\geq 0$ and $\lambda \in T^*M$, and $\lambda = 0 \Leftrightarrow \mathcal L(\lambda) = 0$. As a consequence, $\mathcal L$ is a diffeomorphism between $T^*M_0 := T^*M \setminus \{\text{zero section}\}$ and $TM_0$. One defines the cometric
\begin{equation}\label{}
  \begin{array}{ccc}
    F^* : T^*M \to [0,+\infty), &  & F^* = F \circ \mathcal L.
  \end{array}
\end{equation}
Thus $F^*$ is smooth on $T^*M_0$, $F^*(t\lambda)=tF^*(\lambda) \ \forall t\geq0$, and $F^*(\lambda)=0 \Leftrightarrow \lambda=0$.

On $T^*M$ consider the tautological $1$-form $\alpha_{taut}$ and the canonical symplectic structure $\Omega := d\alpha_{taut}$.

In natural coordinates $(x^1,\dots,p^1,\dots) \simeq \sum_i p^i d{x^i}$ on $T^*M$ associated to a set of coordinates $(x^1,\dots,x^n)$ on $M$, $\alpha_{taut} = \sum_i p^i dx^i$ and $\Omega = \sum_i dp^i \wedge dx^i$. The Hamiltonian $H := \frac{1}{2} (F^*)^2$ induces the Hamiltonian vector field $X_H$ by $-dH = \Omega(X_H,\cdot)$, and one checks that
\begin{equation}\label{spray_geo}
  S := \mathcal L_*X_H
\end{equation}
is a spray. Its flow is the geodesic flow of the Finsler metric $F$, and $S$ is called the geodesic spray.

If $g_{ij}(x,y)$ represents the quadratic form $g_v$~\eqref{quad_form} in natural coordinates, where $v \simeq (x,y)$, and if $g^{ij}$ is the inverse of $g_{ij}$, then the geodesic spray coefficients are
\begin{equation}\label{formula_spray}
\begin{array}{ccc}
  G^i = \frac{1}{2}\gamma^i_{jk}y^jy^k & \text{where} & \gamma^i_{jk} := \frac{1}{2} g^{il} \left( \partial_{x^k}g_{lj} - \partial_{x^l}g_{jk} + \partial_{x^j}g_{kl} \right).
\end{array}
\end{equation}
This is proved by analyzing the Euler-Lagrange equations of the variational problem associated to the integral $\int \frac{1}{2}F^2 dt$. The fundamental difference with Riemannian geometry is that all functions $g^{il}$, $\partial_{x^k}g_{lj}$, etc also depend on the fiber coordinates $y^1,\dots,y^n$, and not only on the $x^1,\dots,x^n$. 

In the context of Finsler metrics, $\Gamma_S$ has a symmetric lift which is more suitable than the Berwald connection. Consider the Cartan tensor: the $(0,3)$-tensor on the bundle $VTM_0$ defined in natural coordinates by
\begin{equation}\label{cartan_tensor}
  A_{ijk} = \frac{1}{4} \frac{\partial^3F^2}{\partial y^i \partial y^j  \partial y^k}.
\end{equation}
Roughly speaking, it governs how $g_v$ varies fiberwise. As explained in~\cite{rad} or in~\cite{bao}, the Chern connection is the symmetric lift of $\Gamma_S$ with coefficients
\begin{equation}\label{chern}
  \Gamma^k_{ij} = \gamma^k_{ij} - g^{is} \left( A_{sjt}\Gamma^t_k - A_{jkt}\Gamma^t_s + A_{kst}\Gamma^t_j \right).
\end{equation}

In the context of Finsler manifolds parallel transport has, as expected, useful metric properties.

\begin{lemma}\label{lemma_metric}
Let $\gamma(t) \subset M$ be a geodesic, and let $V,W$ be vector fields along $\gamma$. Then
\begin{equation}\label{metric_props}
  \frac{d}{dt} g_{\dot\gamma}(V,W) = g_{\dot\gamma}\left( \frac{D_\gamma V}{dt},W \right) + g_{\dot\gamma}\left( V,\frac{D_\gamma W}{dt} \right).
\end{equation}
\end{lemma}

In contrast to the Riemannian case, this formula may not hold when $\gamma$ is not a geodesic. See the appendix for a proof.

A flag pole is a pair $(\Pi,v)$, where $v\in TM_0$ and $\Pi$ is a $2$-plane in $T_{\pi(v)}M$ containing $v$. In the context of Finsler manifolds explained above, the associated flag curvature is
\begin{equation}\label{flagcurv}
  K(\Pi,v) = \frac{g_v(-R^v(w),w)}{g_v(v,v)g_v(w,w) - g_v(v,w)^2}
\end{equation}
where $w\in \Pi$ is any vector such that $\{v,w\}$ is linearly independent. In the Riemannian case, where all the tensors involved only depend on the base point, $K(\Pi,v)$ does not depend on $v \in \Pi$, and is the seccional curvature of $\Pi$.

A vector field $J$ along a geodesic $\gamma$ satisfying $\frac{D^2_\gamma J}{dt^2} - R^{\dot\gamma}(J) = 0$ is called a Jacobi field. This ODE is referred to as the Jacobi equation. The linearization of the geodesic flow can be suitably represented according to the following standard lemma. A proof is found in the appendix.

\begin{lemma}\label{lemma_rep_lf}
Let $\Phi_t$ be the flow of $S$ and set $\zeta(t) = d\Phi_t \cdot \zeta(0)$ where $\zeta(0) \in T_{\dot\gamma(0)}TM$ is fixed. If $J$ is the Jacobi field determined by $\zeta(0) \simeq (\frac{D_{\gamma}J}{dt}(0), J(0))$ under the isomorphism~\eqref{iso_ttm} then
\begin{equation}\label{rep_linflow}
  \zeta(t) \simeq \left( \frac{D_{\gamma}J}{dt}(t), J(t) \right).
\end{equation}
\end{lemma}

The unit sphere bundle $SM = \{v \in TM : F(v)=1 \}$ has contact-type as a hypersurface inside $TM_0$ equipped with the sympletic structure $\Omega_F := (\mathcal L^{-1})^*\Omega$. This is so since $C$ is a Liouville vector field, that is, $L_C \Omega_F = \Omega_F$. Thus $\alpha_F := i_C\Omega_F$ restricts to a contact form on $SM$. The geodesic spray $S$
coincides with the Reeb vector field associated to $\alpha_F|_{SM}$, that is, it satisfies $i_S\alpha_F \equiv1$, $i_S\Omega_F \equiv0$. $SM$ becomes a contact manifold with contact structure $\xi_F = \ker \alpha_F|_{SM} \subset TSM$.

Setting $(\R v)^\bot = \{w\in T_{\pi(v)}M \mid g_v(v,w)=0 \}$ then, since $g_v = \frac{1}{2} \frac{\partial^2F^2}{\partial {y^i}\partial {y^j}} dx^i\otimes dx^j$ and $\alpha_F(v) = \frac{1}{2}\frac{\partial F^2}{\partial y^j}dx^j = \frac{1}{2}\frac{\partial^2 F^2}{\partial y^i\partial y^j}y^jdx^i = g_v(v,d\pi \ \cdot)$ in natural coordinates, one has
\begin{equation}\label{contact_str_split}
  \xi_F|_v \simeq (\R v)^\bot \oplus (\R v)^\bot
\end{equation}
under the isomorphism~\eqref{iso_ttm}. Since the linearized flow $d\Phi_t$ preserves $\xi_F$ we get, in view of~\eqref{contact_str_split} and~\eqref{rep_linflow}, a familiar fact from Riemannian geometry: if $J(0),\frac{D_\gamma J}{dt}(0) \in (\R\dot\gamma(0))^\bot$ then $J(t),\frac{D_\gamma J}{dt}(t) \in (\R\dot\gamma(t))^\bot$, for every $t$.

\subsection{Estimates on the Conley-Zehnder index}\label{proof_hp}

For the sake of completeness we quickly discuss the proof of Theorem~\ref{teo_hp}.
We use the notation established in the last paragraphs.

\subsubsection{A global symplectic trivialization of $\xi_F \subset SS^2$}

Let $F$ be Finsler metric on $S^2$, which is equipped with its orientation induced as a submanifold of $\R^3$. The unit sphere bundle $SS^2$ is equipped with the contact form $\alpha_F$ discussed in~\S~\ref{casefinsler}. For every $v\in SS^2$, $g_v$ is an inner-product on $T_{\pi(v)}S^2$, and there exists a unique vector $v^\bot \in T_{\pi(v)}S^2$ such that $\{v,v^\bot\}$ is a positively oriented $g_v$-orthonormal basis of $T_{\pi(v)}S^2$. Thus $(\R v)^\bot = \R v^\bot$ and we obtain a global trivialization
\begin{equation}\label{global_triv}
  (v,(s,t)) \in SS^2 \times \R^2 \mapsto (sv^\bot,tv^\bot) \in (\R v)^\bot \oplus (\R v)^\bot \simeq \xi_F|_{v}
\end{equation}
of $\xi_F$. The last identification is given by~\eqref{contact_str_split}. One checks easily that this trivialization is symplectic with respect to $d\alpha_F$.

\subsubsection{Estimating the linearized twist}

Let $\gamma(t)$ be a geodesic on $S^2$ with unit speed, and choose $\zeta(0) \in \xi_F|_{\dot\gamma(0)} \subset T_{\dot\gamma(0)}SM$. Let $\Phi_t$ be the Reeb flow of $\alpha_F$ on $SS^2$. Then, setting $\zeta(t) := d\Phi_t \cdot \zeta(0) \in \xi_F|_{\dot\gamma(t)}$, we can use~\eqref{iso_ttm} and identify $\zeta(t) \simeq (\frac{D_\gamma J}{dt}(t),J(t))$ as in~\eqref{rep_linflow}, for some Jacobi field $J$. The vector field $\dot\gamma^\bot$ is parallel, as one can prove by using~\eqref{metric_props} and noting that $g_{\dot\gamma}(\dot\gamma,\dot\gamma^\bot)$ and $g_{\dot\gamma}(\dot\gamma^\bot,\dot\gamma^\bot)$ are constant in $t$. Writing $J = f\dot\gamma^\bot$ we get $\frac{D_\gamma J}{dt} = f'\dot\gamma^\bot$, $\frac{D^2_\gamma J}{dt^2} = f''\dot\gamma^\bot$ and, consequently,
\begin{equation}\label{eq_f}
  f'' = g_{\dot\gamma} \left( \frac{D^2_\gamma J}{dt^2}, \dot\gamma^\bot \right) = g_{\dot\gamma} \left( R^{\dot\gamma}(\dot\gamma^\bot), \dot\gamma^\bot \right)f = -K(T_{\gamma}S^2,\dot\gamma)f.
\end{equation}
Abbreviating $K(T_{\gamma}S^2,\dot\gamma)$ by $K(t)$ we get, after further identifying $$ \zeta(t)\simeq (\frac{D_\gamma J}{dt}(t),J(t)) \simeq (f',f) \stackrel{.}{=} u(t) $$ via~\eqref{global_triv}, an equation
\[
  \dot u = \begin{pmatrix} 0 & -K(t) \\ 1 & 0 \end{pmatrix} u.
\]
Thus, if $\vartheta(t)$ is a smooth lift of the argument of $u(t)$ and the Finsler metric is positively curved then
\begin{equation}\label{estimate_flag}
  \dot\vartheta \geq \delta \text{ for any } \delta \text{ satisfying } 0<\delta <\min\{1,K_{min}\},
\end{equation}
where $K_{min}$ is the infimum among all flag curvatures.

\subsubsection{Rademacher's estimate on the length of a closed geodesic}

Let $l$ denote the infimum among all lengths of closed geodesic loops, and $K_{max}$ be the supremum among all flag curvatures. The following estimate was obtained by Rademacher in~\cite{rad_shortest}, see~\cite{rad} for a detailed account of the subject:
\begin{equation}\label{l_minest}
  0<K_{min} \leq K_{max} \Rightarrow l \geq \frac{\pi\left(1 + \frac{1}{r}\right)}{\sqrt{K_{max}}}
\end{equation}
where $r\geq 1$ is the reversibility. In the Riemannian case $r=1$ and~\eqref{l_minest} is obtained from Klingenberg's estimate on the injectivity radius. 

\begin{lemma}\label{case_s2}
In the particular case of a positively curved Finsler metric on $S^2$, every prime closed geodesic $\gamma$ such that $\dot\gamma$ is contractible in $SS^2$ has length larger than or equal to $2\pi(1+r^{-1})/\sqrt{K_{max}}$.
\end{lemma}

\begin{proof}
$\gamma$ consists of at least two distinct closed loops since, otherwise, $\dot\gamma$ is not contractible in $SS^2$.
\end{proof}

\subsubsection{Estimating the index}

Consider a positively curved Finsler metric $F$ with reversibility $r$ that is strongly $\left({r}/{(r+1)}\right)^2$-pinched. After rescaling we can assume $$ 1=K_{min} \leq K_{max} < \left(\frac{r+1}{r}\right)^2. $$ It follows by Lemma~\ref{case_s2} and  estimate~\eqref{estimate_flag} that $\vartheta$ varies strictly more than $2\pi$ along a closed geodesic $\gamma$ with $\dot \gamma$ contractible in $SS^2$, regardless of the choice of initial condition in $\ker \alpha_F|_{\dot\gamma(0)}$. According to~\eqref{CZ_def} this proves Theorem~\ref{teo_hp}.

\section{Proof of Theorem~\ref{main1}}\label{sec_main1}

We split the arguments in three parts. In \S~\ref{par_topol} we prove Lemma~\ref{lema_topol}. In~\S~\ref{non-exist} we use Theorem~\ref{secao_global_1} to prove (i) in Theorem~\ref{main1}. Finally, in~\S~\ref{examples} we exhibit for any $\epsilon>0$ and $r\geq 1$, examples of Finsler metrics of Randers type on $S^2$ with reversibility $r$ and flag curvatures in $\delta(r) - \epsilon < K \leq 1$ admitting geodesics with one transverse self-intersection.

\subsection{Contact-topological invariants of $k_8$}\label{par_topol}

Let us denote by $\mathcal B_8 \subset \mathcal B$ the set of $C^2$ immersions of $S^1$ into $S^2$ with precisely one transverse self-intersection. Every $c \in \mathcal B_8$ induces an embedded copy of $S^1$ inside $F^{-1}(1) \subset TS^2$ given by $\dot c/F(\dot c)$. Recall the Legendre transform $\mathcal L_F:T^*S^2 \to TS^2$ induced by $F$, which is a continuous fiber-preserving map that restricts to a smooth diffeomorphism $T^*S^2\setminus 0 \simeq TS^2\setminus 0$ and induces a cometric $F^* = F\circ \mathcal L_F$. Then the tautological 1-form $\alpha_{\rm taut}$ restricts to a contact form on $(F^*)^{-1}(1)$. Consequently $\alpha_F := (\mathcal L_F^{-1})^*\alpha_{\rm taut}$ restricts to a contact form on $F^{-1}(1)$, and clearly $\dot c/F(\dot c)$ is positively transverse to the contact structure $\xi_F := \ker\alpha_F$. It is not hard to check that any two $c_0,c_1 \in \mathcal B_8$ are homotopic through curves $c_s \in \mathcal B_8$, $s\in[0,1]$, so that we get a corresponding isotopy $\dot c_s/F(\dot c_s)$ through knots which are transverse to $\xi_F$, thus preserving the knot type and the self-linking number. Consequently, it suffices to exhibit one element $c \in\mathcal B_8$ such that $\dot c/F(\dot c)$ is unknotted and has self-linking number $-1$.

First we prove Lemma~\ref{lema_topol} for the metric $F_0(v) := \sqrt{g_0(v,v)}$, where $g_0$ is the Riemannian metric induced by isometrically embedding $S^2$ in $\R^3$ as the unit sphere. Taking polar coordinates $(\theta,R)$ in $\C$, $\theta \in \R/2\pi\Z$ and $R\geq 0$, consider the embedding $$ \gamma:\C \to S^3 \subset \C^2; \ (\theta,R) \mapsto \frac{1}{\sqrt{1+R^2}}(1,Re^ {i\theta}). $$ Denote by
\[
  \lambda_0 = \frac{1}{4i} \left( \bar z dz - zd\bar z + \bar w dw - wd\bar w \right) = \frac{1}{2} \left( q_0dp_0 - p_0dq_0 + q_1dp_1 - p_1dq_1 \right)
\]
the Liouville form on $\C^2 \simeq \R^4$ with complex coordinates $(z = q_0+ip_0 , w = q_1+ip_1)$. The embedded circles $\gamma_R=\gamma(\cdot,R)$ converge, as $R\to \infty$, to the Hopf fiber $P_0 = \{ \theta \mapsto (0,e^{i\theta}) \mid \theta \in [0,2\pi]\}$ in the $C^1$-topology. All $\gamma_R$ are positively transverse (with respect to $\lambda_0$) to $\ker\lambda_0$. It is well know that $\sl(P_0) = -1$, which implies $\sl(\gamma_R) = -1 \ \forall R>0$. Moreover, $\gamma_{R_0}$ is clearly unknotted since it is the boundary of the embedded disk $\mathcal D_{R_0}$ parametrized by $\gamma|_{\{R\leq R_0\}}$.

Identifying $S^3 \simeq SU(2)$ via $$ (z,w) \simeq \begin{pmatrix} z & w \\ -\bar w & \bar z \end{pmatrix}, \ \ |z|^2 + |w|^2 = 1,  $$ and considering the matrices $$ \begin{array}{ccc} j = \begin{pmatrix} 0 & 1 \\ -1 & 0 \end{pmatrix} & \text{and} & k = \begin{pmatrix} 0 & i \\ i & 0 \end{pmatrix} \end{array} $$ there is a double cover $D:S^3 \to F_0^{-1}(1)$ given by
\[
  A \in SU(2) \simeq S^3 \mapsto (A^{-1}jA,-A^{-1}kA) \in F_0^{-1}(1)
\]
where we see a unit vector $(x,y,t) \in S^2 \subset \R^3$ sitting inside $S^3$ as $(it,u)$ where $u=x+iy$. Here $A^{-1}jA$ represents the base point and $-A^{-1}kA$ represents the tangent vector. We have $D^*\alpha_{F_0} = 4\lambda_0|_{S^3}$ (cf.~\cite{CO,HP}). The factor $4$ appears since a Hopf circle on $S^3$, which has $\lambda_0$-action equal to $\pi$, projects onto the unit velocity vector of a great circle prescribed twice, which has length $4\pi$.

The group of deck transformations of $D$ is precisely $\{id,a\}$, where $a$ is the antipodal map. Thus, for each $R_0>0$, $\Gamma_{R_0} := D \circ \gamma_{R_0}$ is an embedded knot in $F_0^{-1}(1)$ since the curve $\gamma_{R_0}$ does not contain pairs of antipodal points. It is clearly transverse to $\xi_{F_0}$. By the same token, $\widetilde {\mathcal D}_{R_0} := D(\mathcal D_{R_0})$ is an embedded disk with boundary $\Gamma_{R_0}$, proving that $\Gamma_{R_0}$ is unknotted. Moreover, since $D$ is a 1-1 contactomorphism of a neighborhood of $\mathcal D_{R_0}$ in $S^3$ onto a neighborhood of $\widetilde{\mathcal D}_{R_0}$ in $F_0^{-1}(1)$, we get that the self-linking number of $\Gamma_{R_0}$ is also $-1$. We concluded that each $\Gamma_{R_0}$ is unknotted and has self-linking number $-1$ in the contact manifold $(F_0^{-1}(1),\xi_{F_0})$.

It only remains to find $c \in \mathcal B_8$ such that $\Gamma_1$ is transversely isotopic to $\dot c/F_0(\dot c)$. Let $c : [0,2\pi] \to S^2$ be defined by the equation $$ c(\theta) = \pi \circ \Gamma_1(\theta) $$ where $\pi : TS^2 \to S^2$ is the bundle projection. It is clear from the formula
\[
  c(\theta) = \frac{1}{2}(1+\cos 2\theta,\sin 2\theta,2\sin\theta) \in \R^3
\]
that $c\in \mathcal B_8$. Since $\Gamma_1$ is positively transverse to $\xi_{F_0}$, we have  $$ \alpha_{F_0}|_{\Gamma_1(\theta)} \cdot \dot\Gamma_1(\theta) =  g_0(c(\theta))(\dot c(\theta),\Gamma_1(\theta)) > 0 \ \forall \theta. $$ Thus we find a $C^1$ lift $\vartheta(\theta) \in (-\frac{\pi}{2},\frac{\pi}{2})$ for the $g_0$-angle between $\Gamma_1$ and $\dot c/F_0(\dot c)$, and can define a transverse homotopy between $\Gamma_1$ and $\dot c/F_0(\dot c)$ keeping the base points fixed by the formula
\[
  h_s(\theta) = \frac{(1-s)\frac{\dot c(\theta)}{F_0(\dot c(\theta))} + s \Gamma_1(s)}{F_0\left( (1-s)\frac{\dot c(\theta)}{F_0(\dot c(\theta))} + s \Gamma_1(s) \right)} \in T_{c(\theta)}S^2 \cap F_0^{-1}(1).
\]
It remains to show that $\{h_s\}_{s\in[0,1]}$ is a transverse isotopy. The only possibility for self-intersections of the curves $\theta \mapsto h_s(\theta) \in F_0^{-1}(1)$ is at the values $\theta=0$ and $\theta=\pi$ where the curve $c(\theta)$ self-intersects at the point $(1,0,0)$. Looking at the formulas
\[
  \begin{array}{c}
    \dot c(\theta) = (-\sin 2\theta,\cos 2\theta,\cos\theta) \in T_{c(\theta)}S^2 \\
    \Gamma_1(\theta) = (-\frac{1}{2}\sin2\theta,\frac{1}{2}(\cos2\theta-1),\cos\theta) \in T_{c(\theta)}S^2
  \end{array}
\]
we note that both $\dot c(0)$ and $\Gamma_1(0)$ point at north hemisphere, while $\dot c(\pi)$ and $\Gamma_1(\pi)$ point at the south hemisphere. Thus the formula for $h_s$ does not produce self-intersections and, consequently, is a transverse isotopy. Lemma~\ref{lema_topol} is proved for the metric $F_0$.

Now we consider a general Finsler metric $F:TS^2\to [0,+\infty)$. We have associated Legendre transforms $\mathcal L_{F_0},\mathcal L_{F}$ and cometrics $F_0^* = F_0 \circ \mathcal L_{F_0}$, $F^* = F \circ \mathcal L_{F}$. The map $\Psi : (F^*)^{-1}(1) \to (F_0^*)^{-1}(1)$ defined by $\Psi(p) = p/F_0^*(p)$ satisfies $\Psi^*(\alpha_{\rm taut}|_{(F_0^*)^{-1}(1)}) = \frac{1}{F_0^*}\alpha_{\rm taut}|_{(F^*)^{-1}(1)}$, so that it is a contactomorphism. Hence we get a contactomorphism
\[
  \Phi := \mathcal L_{F_0} \circ \Psi \circ \mathcal L_{F}^{-1} : F^{-1}(1) \to F_0^{-1}(1).
\]

Given any immersion $c:S^1 \to S^2$ we construct two transverse embeddings $V_c,W_c : S^1 \to F_0^{-1}(1)$ covering $c$ as follows: $V_c = \dot c/F_0(\dot c)$ and $W_c = \Phi (\dot c/F(\dot c))$.  Since $\Phi$ preserves co-orientations induced by $\alpha_{F}$ and $\alpha_{F_0}$, we have $$ g_0(c(\theta))(V_c(\theta),W_c(\theta)) > 0, \ \forall \theta \in \R/2\pi \Z \simeq S^1. $$ Now let $c_0 \in \mathcal B_8$ be arbitrary. Assume its self-intersection point is $c_0(0) = c_0(\pi)$. We clearly can find a homotopy $c_s:S^1 \to S^2$, $s\in[0,1]$, starting at $c_0$, such that $c_s \in \mathcal B_8$ for every $s\in [0,1)$, and the immersion $c_1$ has a unique self-intersection point at $c_1(0) = c_1(\pi)$ satisfying $\dot c_1(0) = -\dot c_1(\pi)$ (negative self-tangency). This induces corresponding isotopies $V_{c_s},W_{c_s} : S^1 \to F_0^{-1}(1)$ through transverse knots. Now define a homotopy $H_\tau : S^1 \to F_0^{-1}(1)$ satisfying $H_0 = V_{c_1}$ and $H_1 = W_{c_1}$ by
\[
  H_\tau(\theta) = \frac{(1-\tau)V_{c_1}(\theta) + \tau W_{c_1}(\theta)}{F_0((1-\tau)V_{c_1}(\theta) + \tau W_{c_1}(\theta))},
\]
where $\tau \in [0,1]$. The map $H_\tau$ is well-defined since $g_0(c_1(\theta))(V_{c_1}(\theta),W_{c_1}(\theta)) > 0$, for all $\theta$, and hence the denominator above never vanishes. However, there could be some value of $\tau$ where $H_\tau$ is not a knot in $F_0^{-1}(1)$. This would only be the case if $H_\tau$ has self-intersections, which is only possible at the values $\theta=0$ and $\theta=\pi$. Note, however, that this never happens because of the condition $\dot c_1(0) = -\dot c_1(\pi)$. We succeeded in showing that the knots $V_{c_0}$ and $W_{c_0}$ are transversally isotopic in $(F_0^{-1}(1),\xi_{F_0})$. We already showed before that $V_{c_0}$ is unknotted and has self-linking number $-1$. Thus the same is true for the knot $\dot c/F(\dot c) \subset (F^{-1}(1),\xi_F)$.

\subsection{Non-existence of geodesics}\label{non-exist}

Let $F$ be a Finsler metric on $S^2$ satisfying~\eqref{hypothesis1}. Then, as remarked in~\S~\ref{casefinsler}, the pull-back of the tautological $1$-form to $TTS^2$, via the inverse Legendre transform induced by $F$, restricts to a contact form $\lambda$ on $SS^2$, which induces the contact structure $\xi := \ker \lambda$. By Theorem~\ref{teo_hp}, $\lambda$ is dynamically convex.

Assume, by contradiction, that there exists a prime closed geodesic $\gamma$ with unit speed, such that $P = \{ t \mapsto \dot\gamma(t) \} \subset SS^2$ is a closed unknotted Reeb orbit with self-linking number $-1$.

In the case $\lambda$ is non-degenerate Theorem~\ref{secao_global_1} implies that $SS^2$ is homeomorphic to $S^3$, a contradiction. It remains to consider the degenerate case. Denote by $$ \mathcal{O}_{\lambda}= \{f\in C^{\infty}(SS^2,\R^+) \mid f\lambda \mbox{ is non-degenerate}\}. $$ The following lemma, which we state without proof, is an adaptation of Lemma 6.8 in~\cite{convex} to our situation, see also~\cite{props1}. The proof is straightforward.

\begin{lemma}\label{lem_convex}
There exists a sequence $f_k \in \mathcal{O}_{\lambda}$ converging to $1$ in the $C^{\infty}$-topology as $k\to+\infty$, such that each contact form $\lambda_k:=f_k \lambda$ admits $P=(x,T)$ as a closed Reeb orbit.
\end{lemma}

We will prove now that $\lambda_k$ is dynamically convex for all $k$ sufficiently large. We denote by $\phi_t$ the Reeb flow of $\lambda$ and by $\phi^k_t$ the Reeb flow of $\lambda_k$.

Consider the global $d\lambda$-symplectic trivialization $\Psi : \xi \stackrel{\sim}{\rightarrow} SS^2 \times \C$ described in~\eqref{global_triv}. Since $\lambda_k\to\lambda$ in $C^\infty$ we can find $d\lambda_k$-symplectic trivializations $\Psi_k : \xi \stackrel{\sim}{\rightarrow} SS^2 \times \C$ such that $\Psi_k \to \Psi$ in $C^\infty$. Given $p_0\in SS^2$ and $v(0) \in \xi_{p_0}\setminus 0$ arbitrary, the solution $v(t) = d\phi_t \cdot v(0) \in \xi_{\phi_t(p_0)}$ of the linearized $\lambda$-Reeb flow can be represented using the frame $\Psi$ as a curve $v(t) \simeq r(t)e^{i\theta(t)}\subset \mathbb{C}$, where $r(t) >0$ and $\theta(t) \in \R$ is any continuous lift of the argument. There exists $a>0$ such that
\begin{equation}\label{eq_linconv1}
  \dot \theta(t)>a,
\end{equation}
for all $t$, independently of the choice of $p_0$ and $v(0)$. The existence of $a$ follows from estimate~\eqref{estimate_flag}. Moreover, if $p_0$ is a point on a closed contractible $\lambda$-Reeb orbit $P=(\bar x,\bar T)$ then
\begin{equation}\label{eq_linconv2}
  \theta(\bar T)-\theta(0) > 2\pi.
\end{equation}
This is a consequence of~\eqref{CZ_def} and of the dynamical convexity of $\lambda$. In the following, solutions of the linearized $\lambda_k$-Reeb flow will be represented similarly by curves in the complex plane with the use of the frame $\Psi_k$.

Arguing indirectly, suppose there exists a subsequence of $f_k$, again denoted by $f_k$, such that $\lambda_k$ is not dynamically convex. Then there exists a contractible closed $\lambda_k$-Reeb orbit $P_k=(x_k,T_k)$, with $\mu_{CZ}(P_k) \leq 2$, for each $k$. Assume first that $T_k \to +\infty$ as $k\to +\infty$. Since $\lambda_k \to \lambda$ and $\Psi_k\to\Psi$ in the $C^{\infty}$-topology, and $SS^2$ is compact, inequality \eqref{eq_linconv1} holds for any linearized solution $v_k \simeq r_ke^{i\theta_k}$ of the $\lambda_k$-Reeb flow over $P_k$, if $k$ is large. In view of the geometric definition of the Conley-Zehnder index explained in~\S~\ref{CZ-3d}, this implies $$ \mu_{CZ}(P_k) \geq \frac{\theta_k(T_k) - \theta_k(0)}{\pi} -1 > \frac{aT_k}{\pi} -1 \to +\infty $$ as $k\to \infty$, in contradiction with $\mu_{CZ}(P_k)\leq 2$.

Now assume that $T_k = \int_{P_k} \lambda_k$ has a bounded subsequence. By the Arzel\`a-Ascoli theorem we find a converging subsequence, still denoted $T_k$, such that $T_k\to T_0$, $x_k(T_k\cdot)\to x_0(T_0\cdot)$ in $C^{\infty}$ where $P_0=(x_0,T_0)$ is a closed $\lambda$-Reeb orbit.  We also find a solution $v_k=r_ke^{i\theta_k}$ of the linearized flow over $P_k$, with $v_k(0)$ bounded and bounded away from $0\in \C$, satisfying $\theta_k(T_k)-\theta_k(0) \leq 2\pi$ for each $k$. Again using that  $\lambda_k \to \lambda$ and $\Psi_k\to \Psi$ in $C^\infty$, a subsequence of $v_k$ converges (in $C^\infty$) to a solution $v_0=r_0e^{i\theta_0}$ over $P_0$ satisfying $\theta_0(T_0)-\theta_0(0)\leq 2\pi$, in contradiction to~\eqref{eq_linconv2}.

Therefore we end up with a sequence $\lambda_k$ converging to $\lambda$ in the $C^{\infty}$-topology such that, for each $k$ large enough, $\lambda_k$ is a dynamically convex non-degenerate tight contact form on $SS^2$ admitting the unknotted closed orbit $P$ with $sl(P)=-1$. Reasoning as before, Theorem~\ref{secao_global_1} leads to the contradiction $SS^2 \simeq S^3$. The proof of (i) in Theorem~\ref{main1} is complete.

\subsection{Examples of Randers type}\label{examples}

Here we prove (ii) in Theorem~\ref{main1}. A Riemannian metric $a$ and a $1$-form $b$ on a manifold $M$ induce a Randers metric
\begin{equation}\label{randers}
  F(v) = \sqrt{a(v,v)} + b(v), \ v\in TM,
\end{equation}
precisely when $|b|_a <1$ everywhere. These form an interesting and rich family of Finsler geometries, vastly studied in the literature.

\subsubsection{Zermelo navigation}\label{zermelo_nav}

A pair $(h,X)$, where $h$ is a Riemannian metric on $M$ and $X$ is a vector field satisfying $\sqrt{h(X,X)} < 1$, is called a Zermelo navigation data. It induces a Randers-type metric on $T^*M$ by $F^*(\lambda) = \sqrt{h^*(\lambda,\lambda)} + \lambda \cdot X$, where $h^*$ is the dual of $h$. The pull-back $F$ of $F^*$ by the Legendre transform is a Finsler metric on $M$, which is said to solve the associated Zermelo navigation problem. In fact, its geodesic flow parametrizes the movement of a particle on $M$ under the additional influence of a tangential wind, see~\cite{bao} for a detailed discussion.

\begin{remark}\label{legendre_remark}
It is curious that $F$ is of Randers type, that is, Legendre transformation preserves the form of the metric. To see this, consider $(x^1\dots x^n,y^1\dots y^n)$ natural coordinates on $TM$, with dual coordinates $(x^1\dots x^n,p_1\dots p_n)$ on $T^*M$. Then $F^* = \sqrt{h^{rs}p_rp_s} + X^kp_k$ and
\[
  y^i = \frac{1}{2}\frac{\partial(F^*)^2}{\partial p_i} = \left( \frac{h^{ij}p_j}{\sqrt{h^{rs}p_rp_s}} + X^i \right)F^* \Rightarrow p_k = h_{ki}(y^i-FX^i)\frac{\sqrt{h^{rs}p_rp_s}}{F}.
\]
Plugging into the formula for $F = F^*$, and writing $\Delta = p_kX^k$, we get
\[
  F = \frac{F-\Delta}{F} \sqrt{h_{kl}(y^k-FX^k)(y^l-FX^l)} + \Delta
\]
which gives $F = \sqrt{h_{kl}(y^k-FX^k)(y^l-FX^l)}$. Raising to the square and expanding the right side we get a second degree polynomial $\epsilon F^2 + BF + C=0$, with $\epsilon=(1-X_lX^l)$, $B = 2y^kX_k$ and $C=-h_{kl}y^ky^l$.  Here we lowered the indices of $X$ with the metric $h$. Solving we get $F = \sqrt{a_{ij}y^iy^j} + b_ky^k$, with $a_{ij}= \epsilon^{-1}h_{ij} + \epsilon^{-2}X_iX_j$ and $b_k = -\epsilon^{-1}X_k$.
\end{remark}

The behavior of the geodesic flow of $F$ is better understood if we work on $T^*M$ equipped with its canonical symplectic structure $\Omega \simeq dp_i \wedge dx^i$ and Hamiltonian $\frac{1}{2}(F^*)^2$. This discussion is based on~\cite{rad}. We can write
\[
\begin{array}{cccc}
  F^* = H+K & \text{ with } & H = \sqrt{h^{ij}p_ip_j} & K = p_iX^i.
\end{array}
\]
The Hamiltonian vector fields are $X_{F^*} = X_H + X_K$, where $i_{X_H}\Omega = -dH$ and $i_{X_K}\Omega = -dK$. If $R_t$ is the flow of $X$ then the flow of $X_K$ is $\lambda \mapsto (dR_t^{-1})^*\cdot\lambda$. If $X$ is Killing with respect to the metric $h$ then $\{H,K\} = 0$ since $R_t$ are isometries. Thus the flows $\Psi^H_t$ and $\Psi^K_t$ of $X_H$ and $X_K$ respectively, commute and, consequently, $\Psi^{F^*}_t = \Psi^H_t\circ\Psi^K_t$. The geodesic flow is precisely the Hamiltonian flow $\Psi^{(F^*)^2/2}_t$ of $\frac{1}{2}(F^*)^2$. Hence
\begin{equation}\label{geod_flow_formula1}
\begin{array}{cc}
  \Psi^{(F^*)^2/2}_t = \Psi^{F^*}_t = \Psi^H_t\circ\Psi^K_t = \Psi^K_t\circ\Psi^H_t & \text{on the unit sphere bundle}
\end{array}
\end{equation}
since $X_{(F^*)^2/2} = F^* X_{F^*}$. Since $H$ is constant along trajectories of $X_{H^2/2}$ we have $\Psi^H_t \left(\lambda\right) = \Psi^{H^2/2}_{t/H(\lambda)} (\lambda)$, which we can use to finally arrive at
\begin{equation}\label{geodflow_formula2}
  \Psi^{(F^*)^2/2}_t(\lambda) = \Psi^K_t\circ \Psi^{H^2/2}_{t/H(\lambda)} \left(\lambda\right), \ \forall \lambda \in (F^*)^{-1}(1).
\end{equation}

Consider the Legendre transform $\mathcal L(\lambda) \in T_xM$ associated to $\frac{1}{2}(F^*)^2$, where $\lambda = p_idx^i \in T^*_xM$ is some covector. If $F^*(\lambda)=1$ we have
\[
  \mathcal L (\lambda) = \sum_i \frac{1}{2}(\partial_{p_i} (F^*)^2) \partial_{x^i} = \left(\frac{h^{ij}p_j}{\sqrt{h^{rs}p_rp_s}} + X^i\right) \partial_{x^i} = \mathcal L_h \left( \frac{\lambda}{H(\lambda)} \right) + X
\]
where $\mathcal L_h$ is the Legendre transform associated to $\frac{1}{2}H^2$. Thus, setting $\lambda=\mathcal L^{-1}(v)$, we find
\begin{equation}\label{legendre_inverse}
  \mathcal L_h \left( \frac{\lambda}{H (\lambda)} \right) = v - X, \ \forall v\in SM = F^{-1}(1).
\end{equation}

Fix $v\in SM$ and let $c(t)$ be the geodesic of $F$ satisfying $\dot c(0) = v$. Then $c(t)$ is the base point of $\Psi^{(F^*)^2/2}_t(\lambda)$. The base point of $\Psi^{H^2/2}_{t/H(\lambda)} (\lambda )$ is equal to $\gamma_0(t)$, where $\gamma_0$ is the geodesic with respect to the metric $h$ satisfying $\dot\gamma_0(0) = \mathcal L_h(\lambda/H(\lambda)) = v-X$, by~\eqref{legendre_inverse}. This fact and~\eqref{geodflow_formula2} imply
\begin{equation}\label{geodesic_good_form}
  c(t) = R_t (\gamma_0(t)).
\end{equation}
This formula will be used later. The reversibility of $F$ is
\begin{equation}\label{reversibility}
  r = \sup_{x\in M} \frac{1+|X(x)|_h}{1-|X(x)|_h}.
\end{equation}

\subsubsection{Pinched surfaces of revolution}\label{pinched_surf}

Let $(x,y,z)$ be standard Euclidean coordinates in $3$-space. We consider a surface $S$ of revolution around the $z$-axis symmetric under reflection with respect to the $xy$-plane, with the Riemannian metric induced by the ambient Euclidean metric. Denoting $\rho = \sqrt{x^2+y^2}$, then $S$ is determined by a curve in the $\rho z$-plane, which we assume parametrized by a parameter $s$ satisfying $\dot\rho^2+\dot z^2=1$. We make $s=0$ correspond to the equator, where the radius is $R = \rho(0)$. We always assume $R< 1$.

The symmetry condition forces $\rho$ to be an even function of $s$ in its domain $(-L,L)$. If $S$ is a sphere then $L<+\infty$, $\dot\rho \to -1$ as $\rho \to \pm L$, the length of a meridian is $2L$ and $S$ intersects the $z$-axis in two poles. In the complement of the poles we have obvious coordinates $(s,\theta) \in (-L,L) \times \R/2\pi \Z$.

The Gaussian curvature is $K = -\ddot\rho/\rho$, and we denote by $K_{max}$ and $K_{min}$ its maximum and minimum, respectively. If $K$ is everywhere positive then the maximal radius is attained at the equator, and the maximal height is attained at the poles. We wish to construct $S$ in a way that $K \equiv K_{min} = 1$ holds at the equator. Then $\ddot \rho(0) = -R$ and, moreover, $K_{max} \geq 1/R^2$ is a necessary condition. In fact, using $\dot\rho(0)=0$, we compute
\begin{equation*}\label{}
  \begin{aligned}
    1-R^2 &= \int_0^L (\rho^2+\dot\rho^2)' ds = \int_0^L (K-1)(-2\rho\dot\rho)ds \\
    &\leq (K_{max}-1) \int_0^L(-2\rho\dot\rho)ds = (K_{max}-1)R^2.
  \end{aligned}
\end{equation*}

Now we claim that if $K_{max} > 1/R^2$ then $S$ with all the above properties exists. To see that, consider a smooth function $g:[0,R^2] \to [0,+\infty)$ satisfying $g(0)=1$, $g(R^2) = 0$, $g' \in [-K_{max},-1]$, $g'\equiv -K_{max}$ near $0$, $g' \equiv -1$ near $R^2$ and $g''\geq 0$. Note that $g(x) = R^2-x$ for $x \sim R^2$. Let $\rho(s)$ be the unique solution of
\begin{equation}\label{ode}
  \dot\rho = -\sqrt{g(\rho^2)}, \ \ \rho(0)=R
\end{equation}
for $s\geq 0$, which coincides\footnote{Here it should be noted that~\eqref{ode} does not have a unique solution, as one can see by considering the constant function $R$.} with $R\cos s$ for $s$ small. Then $\dot\rho(0)=0$ and there exists $L>0$ such that $\rho(s) \to 0^+$ and $\dot\rho(s) \to -1^+$ as $s\to L^-$. This solution $\rho$ determines a $C^1$-embedded disk in the half-space $z\geq 0$, which can be reflected to provide the required $C^1$ sphere of revolution $S$. One can check that $S$ is smooth, that the maximal value of the Gaussian curvature is $K_{max}$ (attained around poles), and that the minimal value $K_{min} = 1$ (attained around the equator).

Note that $g(x) \geq R^2-x \ \forall x\in [0,R^2]$. This implies
\begin{equation}\label{est_above}
  \rho(s) \leq R\cos s, \ \forall s\in [0,L].
\end{equation}

\subsubsection{Estimating a return time}\label{return_time}

Fix any $R\in (0,1)$ and let $g:[0,R^2] \to [0,+\infty)$ be a function as in \S\S~\ref{pinched_surf} with $K_{max}>R^{-2}$. Consider the associated unique solution $\rho(s):[0,L]\to \R$ of~\eqref{ode} which equals $R\cos s$ for small values of $s$. We claim that, for any $b>1$, it is possible to make $2L < b\pi R$ by taking $K_{max}$ close enough to $R^{-2}$.

To prove this, let $h:[0,R^2] \to [0,+\infty)$ be the continuous function defined by $h(x)=1-K_{max}x$ if $x\in [0,x_*]$ and $h(x)=R^2-x$ if $x\in[x_*,R^2]$, where $x_*=(1-R^2)/(K_{max}-1)\in(0,R^2)$. It is imediate that $h(x)=g(x)$ for all $x$ in a neighborhood of $\{0,R^2\}$. Since $g''\geq0$, we have $g\geq h$.

Let $\xi(s)$ be unique solution of $\dot \xi = -\sqrt{h(\xi^2)}$ with initial condition $\xi(0)=R$ coinciding with $R\cos s$ when $s>0$ is small. We want to estimate the first $s_*>0$ such that $\xi(s_*) =0$. Observe that $\xi$ satisfies $\ddot \xi = -K_{max} \xi$ if $0\leq \xi \leq \sqrt{x_*}$ and satisfies $\ddot \xi = - \xi$ if $\sqrt{x_*} \leq \xi \leq R$. This implies that $s_*=s_1+s_2$ where $s_1$ is such that $\xi(s_1)=\sqrt{x_*}$.

Now we prove that if we choose $K_{max}$ sufficiently close to $1/R^2$ then $s_*$ is close to $\frac{\pi R}{2}$. Observe that $\xi(s) = R\cos s$ for $0\leq s\leq s_1$. Thus if $K_{max} \to R^{-2}$, $x_* \to R^2$ and, therefore,  $s_1 \to 0$. For $s\geq s_1$, $\xi(s)$ is a solution of $\ddot\xi = -K_{max}\xi$ satisfying $\xi(s_1) = \sqrt{x_*}$ and $\dot\xi(s_1)<0$. Thus the time $s_2$ that it takes to reach zero is smaller than the time $R\cos(\sqrt{K_{max}}s)$ takes to decay from $R$ to $0$, which is $\frac{\pi}{2\sqrt{K_{max}}}$. Consequently $$ s_2 < \frac{\pi}{2\sqrt{K_{max}}} \to \frac{\pi R}{2} $$ as $K_{max} \to 1/R^2$. Thus $K_{max} \sim R^{-2}$ implies $s_* = s_1+s_2 < b\pi R/2$.

To estimate the length of the meridian observe that $\rho(s) \leq \xi(s)$ for all $s$ since $g\geq h$ on $[0,R^2]$. Thus the length $2L$ of the meridian is at most $2s_*$ which is smaller than $b\pi R$ for $K_{max}$ close enough to $R^{-2}$.

\subsubsection{Introducing the wind and completing the proof}

So far we have not fixed any of the data explicit in the statement of the Theorem~\ref{main1}.

Let $r\geq 1$ be given. Consider $\epsilon>0$ small and numbers $R$, $K_{max}$ satisfying
\begin{equation}\label{pinch_R}
  \frac{r}{r+1}-\epsilon < R < \frac{r}{r+1}, \ \ \left( \frac{r}{r+1} \right)^{-2} < \frac{1}{R^2} < K_{max} < \left( \frac{r}{r+1}-\epsilon \right)^{-2}.
\end{equation}
Following the construction in \S\S~\ref{pinched_surf}, we can find a smooth surface of revolution $(S,h)$ with Gaussian curvature taking values in $[1,K_{max}]$, and with an equator of radius $R$. We can arrange so that the curvature equals $K_{max}$ near the poles, and equals $1$ near the equator.

By the discussion of \S\S~\ref{return_time} we can assume, after making $\epsilon$ small enough, that the length $T$ of a meridian satisfies
\begin{equation}\label{estimate_time}
  T \leq \left(\frac{r+1}{r}\right)\pi R < \pi.
\end{equation}
We can also assume that $(\frac{r-1}{r+1})/(\frac{r}{r+1}-\epsilon) \sim \frac{r-1}{r} <1$ by making $\epsilon$ even smaller. Take $\eta\geq0$ so that $\eta R = \frac{r-1}{r+1}$. Similarly to~\cite{katok}, consider the vector field
\begin{equation}\label{wind}
  X = \eta \frac{\partial}{\partial \theta}
\end{equation}
and let $F$ be the Randers metric on $S$ induced by the navigation data $(h,X)$, as explained in \S~\ref{zermelo_nav}. By~\eqref{reversibility} $F$ has reversibility $r$. Later we shall need to note that
\begin{equation}\label{est_T_eta}
  T\eta \leq \pi \left( \frac{r-1}{r+1}\right ) \left( \frac{r}{r+1}-\epsilon \right)^{-1} < \pi
\end{equation}

Crucial to our analysis is the fact that all flag curvatures of $F$ are independent of the chosen flagpole and coincide with the Gaussian curvatures of $h$, see~\cite{bao}.

\begin{lemma}
Let $x$ be a point in the equator and let $v \in T_xS$ satisfy $F(v)=1$ and $h(v-X,X)>0 \ (\geq0)$. Then the geodesic $c(t)$ with respect to $F$ with initial condition $\dot c(0) = v$ satisfies $h(\dot c(t),X\circ c(t))>0 \ (\geq 0)$, $\forall t\geq0$.
\end{lemma}

\begin{proof}
Clearly $X$ is Killing for $h$. According to~\eqref{geodesic_good_form}, $c(t) = R_t \circ \gamma_0(t)$, where $R_t$ is the flow of $X$ and $\gamma_0$ is a geodesic with respect to $h$ with $\dot\gamma_0(0)=v-X$. Thus, in view of the Clairaut integral for surfaces of revolution, $h(\dot\gamma_0,X\circ\gamma_0) >0 \ (\geq 0)$ for $t\geq0$. We can estimate
\[
  \begin{aligned}
    h(\dot c,X\circ c) &= h(X\circ c,X\circ c) + h(dR_t \cdot \dot\gamma_0,X \circ c) \\
    & = |X\circ c|_h^2 + h(dR_t \cdot \dot\gamma_0, dR_t \cdot X \circ \gamma_0) \\
    & = |X\circ c|_h^2 + h(\dot\gamma_0, X \circ \gamma_0) >0 \ (\geq 0).
  \end{aligned}
\]
\end{proof}

Fix a point $x_0$ in the equator and let $0<\phi_0\leq\pi/2$ be determined as follows: the unique vector $v_{\phi_0} \in T_{x_0}S$ pointing to the northern hemisphere, satisfying $F(v_{\phi_0})=1$ and $h(v_{\phi_0} - X,X)=0$ makes $h$-angle $\phi_0$ with $X$. If for every $\phi \in [0,\phi_0]$ we denote by $v_\phi$ the unique vector not pointing south, satisfying $F(v_\phi)=1$ and making $h$-angle $\phi$ with $X$, then $\phi<\phi_0 \Rightarrow h(v_\phi-X,X)>0$. Analogously, we write $c_\phi$ for the geodesic of $F$ satisfying $\dot c_\phi(0)=v_\phi$, and $\theta_\phi(t)$ for the unique lift of the function $\theta\circ c_\phi(t)$ to the universal covering $\R$ satisfying $\theta_\phi(0)=0, \ \forall \phi\in [0,\phi_0)$. The lemma above provides the estimate $\dot\theta_\phi>0$.

By~\eqref{geodesic_good_form}, $c_{\phi_0}(t) = R_t \circ \gamma_+(t)$ where $\gamma_+$ is a geodesic of $h$ heading north that leaves $x_0$ $h$-perpendicularly to the equator, and $R_t$ is the flow of $X$. Thus $c_{\phi_0}$ passes through the north pole, and that is why $\theta_\phi(t)$ is defined only for $\phi <\phi_0$. Moreover, $1 = F(v_{\phi_0}) = |v_{\phi_0}-X|_h = |\dot\gamma_+(0)|_h$, see Remark~\ref{legendre_remark}.

Let $T_\phi>0$ denote the first return time of the geodesic $c_\phi$ to the equator, and $P(\phi)$ be the point of return, which are well-defined smooth functions of $\phi \in (0,\phi_0]$. This is so since, by uniqueness of solutions of ODEs, the first hit of any geodesic $c_\phi$, with $\phi\in (0,\phi_0]$, with the equator is transverse. We have that $T_{\phi_0}$ is equal to the $h$-length of the meridian $T$. To see that one needs to make use of the identity $F(v) = |v-X|_h$. For $\phi < \phi_0$ we have $$ P(\phi) = (R\cos \theta_\phi(T_\phi), R\sin \theta_\phi(T_\phi),0). $$

Clearly the formula $c_{\phi_0}(t) = R_t \circ \gamma_+(t)$ implies $$ P(\phi_0) = (R\cos (\pi+T_{\phi_0}\eta), R\sin (\pi+T_{\phi_0}\eta),0). $$  Thus $\text{dist}(\theta_\phi(T_\phi),\{\pi+T_{\phi_0}\eta+ \{0,2\pi,4\pi,\dots\} \}) \to 0 $ as $\phi\to\phi_0^-$. We used that $\dot\theta_\phi>0$ for $\phi<\phi_0$. The curves $c_\phi$ converge in $C^1_{loc}$ to $c_{\phi_0}$ as $\phi\to\phi_0^-$, and $c_{\phi_0}$ does not self-intersect before hitting the equator since $T_{\phi_0}\eta < \pi$ by~\eqref{est_T_eta}. Thus $c_\phi$ does not have a self-intersection before first hitting the equator when $\phi$ is close to $\phi_0$. This proves that $\theta_\phi(T_\phi) \to \pi+T_{\phi_0}\eta$ as $\phi\to\phi_0^-$.

The geodesic flow is the Reeb flow in the unit sphere bundle $SS^2 = F^{-1}(1)$ equipped with the contact form $\alpha_F$ discussed in \S~\ref{geomsetup}. The Jacobi vector field $J(t) = \left.\partial_\phi\right|_{\phi=0} c_\phi(t)$ along $c_0(t)$ satisfies $g_{\dot c_0}(\dot c_0,J(t))\equiv0$ because it comes from a vertical variation and, as such, lies in the contact structure $\xi_F = \ker \alpha_F$. Thus $J(t)=f(t)\dot c_0^\bot(t)$, where $f$ satisfies~\eqref{eq_f}. Here $g_v$, $v\not=0$, is the positive inner-product~\eqref{quad_form} on $T_{\pi(v)}S^2$, and $\dot c_0^\bot$ is the unique vector such that $\{\dot c_0,\dot c_0^\bot\}$ is a positively oriented $g_{\dot c_0}$-orthonormal basis of $T_{c_0}S^2$. Since the flag curvatures along the equator are constant equal to $1$, $f(t)$ is a (positive) multiple of $\sin t$ and, consequently, the first zero of $J$ appears at time $\pi$. Thus $$ P(\phi) \to c_0(\pi) = (R\cos(\pi R^{-1}+\pi\eta) , R\sin(\pi R^{-1}+\pi\eta) ,0) $$ and $\theta_\phi(T_\phi) \to \pi R^{-1}+\pi\eta$ as $\phi\to 0^+$. Summarizing, we proved
\begin{equation}\label{est11}
  \lim_{\phi\to\phi_0^-} \theta_\phi(T_\phi) = \pi + T_{\phi_0}\eta < 2\pi \ \ \ \text{(see~\eqref{est_T_eta})}
\end{equation}
and
\begin{equation}\label{est21}
  \lim_{\phi\to0^+} \theta_\phi(T_\phi) = \pi R^{-1}+ \eta \pi = R^{-1}\pi (1+\eta R) > \left( \frac{r}{r+1} \right)^{-1} \pi \frac{2r}{r+1} = 2\pi.
\end{equation}
By continuity of $\phi\mapsto \theta_\phi(T_\phi)$, there exists $\phi^* \in (0,\phi_0)$ such that $\theta_{\phi^*}(T_{\phi^*})=2\pi$, and $c_{\phi^*}$ first returns to the equator exactly at $x_0$. Moreover, $c_{\phi^*}$ does not self-intersects before hitting the equator since, otherwise, there would be some $t< T_{\phi^*}$ such that $\dot\theta_{\phi^*}(t) \leq 0$, a contradiction.

By the symmetry of $F$ under the reflection $Q$ with respect to the $xy$-plane, $$ \dot c_{\phi^*}(T_{\phi^*}) = dQ \cdot \dot c_{\phi^*}(0). $$ Here one has to make use of the Clairaut integral for the underlying Riemannian metric $h$. Thus, $c_\phi^*:[0,2T_{\phi^*}]\to S^2$ is a smooth closed geodesic with precisely one transverse self-intersection. The flag curvatures lie between $1$ and $K_{max}$, where $K_{max}$ satisfies~\eqref{pinch_R}. We can normalize the curvature, after dilating the Finsler metric, in order to complete the proof of Theorem~\ref{main1}.

\begin{remark}
Given $r\geq 1$ and $0<\delta <(r/(r+1))^2$ we can construct a surface of revolution $S_{\delta,r}$ as above to find an example of a $\delta$-pinched Finsler metric on $S^2$ with reversibility equal to $r$, which is not dyna\-mically convex. In fact, the double cover of the equator of $S_{\delta,r}$ corresponds to a contractible closed geodesic on its unit tangent bundle and has Conley-Zender index equal to $1$. Therefore, the pinching condition on the flag curvatures given by Harris-Paternain in Theorem~\ref{teo_hp} that ensure dynamical convexity for $\alpha_F$ is sharp.
\end{remark}

\section{Proof of Theorem~\ref{main2}}

To prove the first assertion, observe that we can assume by contradiction the existence of convex hypersurfaces $S_k\subset \R^{2m}$, $k\geq 1$, converging to $S^{2m-1}$ in $\textbf{Conv}(2m)$ as $k \to \infty$, such that the Hamiltonian flow on $S_k$ admits a closed orbit $P_k$ with $\A(P_k)=\int_{P_k} \lambda_0 \to T_0\in \R^{+} \backslash \{\pi\}$. These corresponds to the existence of functions $f_k:S^{2m-1}\to (0,+\infty)$ converging to $1$ in the $C^2$-topology such that the Reeb flow $\Phi^k:\R \times S^{2m-1} \to S^{2m-1}$ associated to the contact form $\lambda_k:=f_k \lambda_0|_{S^{2m-1}}$ admits a closed orbit, also denoted by $P_k$, with prime period $T_k$ and satisfying $\A(P_k)=\int_{P_k} \lambda_k =T_k \to T_0\in \R^{+} \backslash \{\pi\}$.  The Reeb vector fields $R_0$, $R_k$, $k\geq 1$, associated to $\lambda_0|_{S^{2m-1}}$, $\lambda_k$, $k\geq 1$, satisfy $R_k \to R_0$, as $k\to \infty$, in the $C^1$-topology. Since all the orbits of the Reeb flow $\Phi^0$ associated to $\lambda_0|_{S^{2m-1}}$ are closed with prime period $\pi$, and since $T_k \to T_0$, we conclude from Arzel\`a-Ascoli theorem the existence of a simple closed orbit $P_0\subset S^{2m-1}$ of $\Phi^0$ and a subsequence of $P_k$, again denoted by $P_k$, such that $P_k$ converges to a $k_0$ cover of $P_0$, for some integer $k_0>1$ satisfying $T_0=k_0 \pi$.  The following proposition due to Bangert~\cite{bangert} is crucial for completing the proof.

\begin{prop}[Bangert]\label{propBangert}
Let $\Phi:\R \times M \to M$, with $M$ closed, be a $C^1$ flow and $p\in M$ be a periodic point of $\Phi$ with prime period $T$. Then for every $\epsilon>0$ there exists a neighborhood $\U$ of $\Phi$ in the weak $C^1$ topology in $C^1(\R \times M, M)$ and a neighborhood $U$ of $p$ in $M$ such that the following holds: if a flow $\hat \Phi \in \U$ has a periodic point $\hat p\in U$ with prime period $\hat T$ then either $\hat T  >\epsilon^{-1}$ or there exists an integer $\hat k>0$ such that $|\hat T - \hat k T|< \epsilon$ and the eigenvalues of the linear map \begin{equation}\label{eq_linear} D_2\Phi(T,p):T_p M \to T_p M \end{equation} which are $\hat k$'th roots of unit generate all the $\hat k$'th roots of unity.
\end{prop}

Applying Proposition \ref{propBangert} to our situation we conclude that
\begin{equation}\label{eq_linear2}
  D_2\Phi^0(\pi,p_0)=Id_{T_{p_0}S^{2m-1}}:T_{p_0}S^{2m-1} \to T_{p_0}S^{2m-1},p_0\in P_0,
\end{equation}
admits an eigenvalue $\lambda_0$ which generates all the $k_0$'th roots of unity. But this is a contradiction since $k_0>1$ and all eigenvalues of \eqref{eq_linear2} are equal to $1$.

To prove the statement made in the case $2m=4$ about the linking numbers of short and long orbits we proceed indirectly and assume, by contradiction, the existence of $S_n \to S^3$ in $\textbf{Conv}(4)$, $P_n,P_n' \in \P(S_n)$ satisfying $\A(P_n) \to \pi$, $\A(P_n') \to +\infty$ and $\sup_n \link(P_n',P_n) =k_0 < \infty$.

We may view $P_n$ and $P_n'$ as closed Reeb orbits in $S^3$ of contact forms $f_n\lambda_0|_{S^3}$ with $f_n \to 1$ in the $C^2$-topology. Let $R_0$ and $R_n$ be the Reeb vector fields associated to $\lambda_0|_{S^3}$ and $f_n\lambda_0|_{S^3}$, respectively. Then $R_n \to R_0$ in the $C^1$-topology. We denote by $\Phi^n_t$ and $\Phi^0_t$ the flows of $R_n$, $R_0$ respectively. We can also assume, in view of the Arzel\`a-Ascoli theorem, that $P_n\to L_0$ in $C^1$ for some Hopf fiber $L_0 \subset S^3$.

The Hopf fiber $L_0$ corresponds to a $\pi$-periodic orbit of the flow $\Phi^0$. Identifying $\R^4 \simeq \C^2$, there is no loss of generality if we assume $L_0 = \{(e^{i\vartheta},0) \mid \vartheta \in \R\}$. Denote by $E:\C^2 \setminus \{0\} \to S^3$ the projection $E(x) = x/|x|$ and by $\Pi_\tau$ the plane $\{e^{i2\pi \tau}\}\times \C$. Then each $D_\tau := \overline {E(\Pi_\tau)}$ is an embedded disk, transverse to the vector field $R_0$ in $\interior{D_\tau}$, satisfying $\partial D_\tau = L_1$, where $L_1$ is the Hopf fiber $L_1 = \{(0,e^{i\vartheta}) \mid \vartheta \in \R\}$. Moreover, $D_\tau$ is a global surface of section for $\Phi^0$ and the (first) return map to $\interior{D_\tau}$ is precisely the identity.

This open book decomposition induces a diffeomorphism
\begin{equation}\label{}
  \begin{aligned}
    \Psi : S^3 \setminus L_1 &\to \R/\pi \Z \times \C \\
    (z,w) &\mapsto \left( \frac{\arg z}{2}, \frac{w}{\sqrt{1-|w|^2}} \right)
  \end{aligned}
\end{equation}
satisfying $\Psi(\interior{D_\tau}) = \{\tau\} \times \C$. The flow $\hat \Phi^0_t := \Psi_*\Phi^0_t = \Psi \circ \Phi^0_t \circ \Psi^{-1}$ is obtained by integrating the vector field $\hat R_0 := \Psi_*R_0$ and is given by
\begin{equation}\label{round_flow}
  \hat \Phi^0_t(\tau_0,\zeta) = (\tau_0+t, e^{i2t}\zeta).
\end{equation}
Note that $L_0$ is mapped precisely onto $\R/\pi\Z \times 0$ since $\Psi(e^{i2\vartheta},0) = (\vartheta,0)$.

Since $P_n \to L_0$ in $C^1$, we find diffeomorphisms $F_n$ of $\R/\pi\Z \times \C$ satisfying: $\supp (F_n) \subset \R/\pi\Z\times B_1(0)$, $F_n \to id$ in $C^1$ and $F_n \circ \Psi(P_n) = \R/\pi\Z \times 0$. This follows from an application of Lemma~\ref{lemmaP_n}.

\begin{lemma}\label{lemmaP_n}
Let $M$ be a $C^k$-manifold and $Z\hookrightarrow M$ be a closed $C^k$-submanifold. Suppose $U$ is any open neighborhood of $Z$ and $Z_n$ are $C^k$-submanifolds converging to $Z$ in the $C^k$-topology as $n\to\infty$. Then there exist diffeomorphisms $\varphi_n$ satisfying $\varphi_n(Z_n)=Z$, $\supp(\varphi_n) \subset U$ and $\varphi_n \to id_M$ in the $C^k$-topology.
\end{lemma}

\begin{proof}
By considering a tubular neighborhood of $Z$ in $M$ one sees that there is no loss of generality if we assume $M$ is $C^k$-vector bundle, $Z$ is the zero section, $U$ is a neighborhood of the zero section and the $Z_n$ are graphs of sections $s_n$ converging to the zero section in the $C^k$-topology. Let $f:M \to[0,1]$ be a fixed smooth function with support compactly contained in $U$ that is identically equal to $1$ in a neighborhood of the zero section containing all $Z_n$. If $\pi$ is the projection of $M$ onto its base, we define the diffeomorphism $\varphi_n : M \to M$ by $\varphi_n(v) := v- f(v)s_n(\pi(v))$. It is easy to check that $\varphi_n$ satisfies all requirements, when $n$ is sufficiently large.
\end{proof}

If we set $\Psi_n := F_n \circ \Psi$, $\hat R_n := (\Psi_n)_* R_n$ and denote by $\hat\Phi^n$ the flow of $\hat R_n$ then the maximal domain of definition of $\hat\Phi^n$ is an exhausting sequence of open subsets of $\R\times \R/\pi\Z \times \C$ and $\hat\Phi^n \to \hat\Phi^0$ in the $C^1$-topology on compact sets of $\R\times \R/\pi\Z \times \C$.

Consider $D_0 := \overline {\Psi^{-1}(\R/\pi\Z \times [0,+\infty))}$. This set is a smooth embedded disk satisfying $\partial D_0=L_0$ which is transverse to $R_0$ at $\interior{D_0} = D_0 \setminus L_0$. In fact, $\interior{D_0}$ coincides with $E(\C\times \{1\})$. Now we consider disks
\begin{equation}\label{}
  D_n := \overline {\Psi_n^{-1}(\R/\pi\Z \times [0,+\infty))}.
\end{equation}
Since the support of $F_n$ is contained in $\R/\pi\Z \times B_1(0)$ the disk $D_n$ coincides with $D_0$ on $S^3 \setminus \Psi^{-1}(\R/\pi\Z \times B_1(0))$. Moreover, $\partial D_n = P_n$.

We claim $D_n$ is transverse to $R_n$ at the points of $D_n \setminus P_n$.
Using Taylor's formula and comparing with $\hat R_0$ we obtain
\begin{equation}\label{}
  \begin{aligned}
    \hat R_n(\tau,z) &= \hat R_n(\tau,0) + D_2\hat R_0(\tau,0) \cdot z \\
    &+ \left[ D_2 \hat R_n(\tau,0)- D_2 \hat R_0(\tau,0) + \int_0^1 D_2 \hat R_n(\tau,\lambda z)- D_2 \hat R_n(\tau,0) d\lambda \right] \cdot z \\
    &= \hat R_n(\tau,0) + D_2\hat R_0(\tau,0) \cdot z + \epsilon_n(\tau,z) \cdot z
  \end{aligned}
\end{equation}
with
\[
  \sup \{ |\epsilon_n(\tau,z)| : (\tau,z) \in \R/\pi\Z \times B_1(0) \} \to 0
\]
as $n \to \infty$. In these coordinates the vector $(0,0,1)$ is normal to the strip $\R/\pi\Z \times [0,+\infty)$. Note that
\begin{equation}\label{DR_n}
  D_2 \hat R_0 (\tau,0) = \begin{bmatrix} 0 & 0 \\ 0 & -2 \\ 2 & 0 \end{bmatrix}
\end{equation}
so for every $\eta>0$ we find $r_0>0$ and $n_0>0$ such that
\begin{equation}\label{}
  \left< \hat R_n(\tau,r) , (0,0,1) \right> \geq (2-\eta)r
\end{equation}
for every $n\geq n_0$ and $r<r_0$. This shows that $\hat R_n$ is transverse to $\R/\pi\Z \times (0,r_0)$ if $n$ is large since we can choose $\eta <2$, which amounts to say that $\Psi_n^{-1}(\R/\pi\Z \times (0,r_0))$ is transverse to $R_n$. Now observe that $\overline {\Psi_n^{-1}(\R/\pi\Z \times [r_0,+\infty))}$ is converging in the $C^1$-topology to the disk $\overline {\Psi^{-1}(\R/\pi\Z \times [r_0,+\infty))}$, which is transverse to $R_0$. So $\overline {\Psi_n^{-1}(\R/\pi\Z \times [r_0,+\infty))}$ will also be transverse to $R_n$ if $n$ is large, proving our claim.

\begin{remark}
Fix $n$ and suppose $\gamma:[a,b] \to S^3 \setminus (P_n\cup L_1)$ is a closed curve. Define $s\in[a,b] \mapsto \zeta(s) \in \C\setminus\{0\}$ by $\Psi_n \circ \gamma(s) = (*,\zeta(s))$. Then it is not hard to check that
\begin{equation}\label{link_formula}
  \link(\gamma,P_n) = \frac{\vartheta(b)-\vartheta(a)}{2\pi}
\end{equation}
where $\vartheta:[a,b] \to \R$ is any continuous lift of the argument of $\zeta(s)$.
\end{remark}

We split the remaining arguments in two cases. \\

\noindent {\bf Case 1:} $\forall n \ \exists y_n \in P_n'$ satisfying $\dist(y_n,L_1) \to 0$. \\

Since $L_1$ is a closed orbit of $R_0$, we can assume the existence of $z_1\in L_1 \cap D_0$ and $y_n \in P_n'$ such that $y_n \to z_1$ as $n \to \infty$.  Let $V\subset D_0$ be a small neighborhood of $z_1$. Note that by the construction above, $V\subset D_n$ for all $n$ large. Now since $R_n \to R_0$ and $\Phi^0_{\pi}(z_1)=z_1$, given any integer $k>0$ and any real number $\epsilon>0$ we can find $U_k\subset V$ neighborhood of $z_1$, $n_0(k,\epsilon)>0$, both depending on $k$ and $\epsilon$, such that for all $n>n_0(k,\epsilon)$ and $z\in U_{k,\epsilon}$, the solution $\{\Phi^n_t(z)$, $t\in(0,(k+1/2)\pi]\}$, intersects $V$ transversely and positively at least $k$ times, and these intersections  correspond to points $t\in [j\pi-\epsilon,j\pi+\epsilon]$ for each $j\in\{1,\ldots,k\}$. Now since $\A(P_n') \to \infty$ as $n\to \infty$, it follows that $\link(P_n ,P_n')=\#\{P_n' \cap D_n\} \to \infty$ as $n\to \infty$ which is a contradiction. In this last assertion we strongly used that each $D_n$ is a (positively) transverse disk to $R_n$, when $n$ is large, so $P_n'$ never intersects $D_n$ negatively. \\

\noindent {\bf Case 2:} $\inf \{\dist(y,L_1) \mid y\in P_n' \} \geq \delta > 0$. \\

By our hypotheses $\exists \rho>0$ such that $\Psi(P_n') \subset \R/\pi\Z \times B_\rho(0)$ when $n$ is large enough. Moreover, $\hat P_n' := \Psi_n(P_n')$ is a $T_n'$-periodic orbit of the flow $\hat\Phi^n_t$ completely contained in $\R/\pi\Z \times B_\rho(0)$, and $T_n' \to \infty$.

For each $n$ fix a point $y_n \in \hat P_n' \cap (0 \times B_\rho(0))$. Define $\gamma_n(t) := \hat\Phi^n_t(y_n) = (\tau_n(t),\zeta_n(t)) \in \R/\pi\Z \times B_\rho(0)$. If $\tilde\tau_n(t)$ is a continuous lift of $\tau_n(t)$ and $\vartheta_n(t)$ is a continuous lift of the argument of $\zeta_n(t)$, then we consider
\[
\begin{array}{ccc}
  N_n := \frac{\tilde\tau_n(T_n')-\tilde\tau_n(0)}{\pi} \in \Z^+ & \text{and} & l_n := \frac{\vartheta_n(T_n')-\vartheta_n(0)}{2\pi} \in \Z.
\end{array}
\]
By~\eqref{link_formula} we have $l_n = \link(P_n',P_n)$.

The (first) return time with respect to the flow $\hat\Phi^n_t$ for points of $0\times B_\rho(0)$ to return to $0 \times \C$ is a function converging uniformly to the constant $\pi$ on $0\times B_\rho(0)$ since $\hat\Phi^n \to \hat\Phi^0$ as above. Since $\hat R_n \to \hat R_0$ in $C^1_{loc}$, $\hat R_n$ is transverse to the disks $\tau \times B_\rho(0)$, $\tau \in \R/\pi\Z$, when $n$ is large enough. Consequently we can divide, for each $n$ large enough, the interval $[0,T_n']$ in precisely $N_n$ intervals $\{I^n_1,\dots,I^n_{N_n}\}$ of lengths uniformly close to $\pi$ corresponding to points where $P_n'$ intersects $0 \times B_\rho(0)$. This implies that $N_n \to \infty$ as $n\to\infty$.

We need to estimate $\dot \vartheta_n$. Let us write $\hat R_n = (a_n,Y_n=(u_n,v_n))$. If $\zeta_n(t) = x_n(t) + i y_n(t)$ then
\[
  \begin{aligned}
    \dot \vartheta_n &= \frac{-y_n\dot x_n + x_n\dot y_n}{x_n^2+y_n^2} = \frac{\left< \zeta_n , -J_0 \cdot \dot\zeta_n \right>}{|\zeta_n|^2} \\
    &= \frac{\left< \zeta_n , -J_0 \cdot (Y_n\circ \gamma_n) \right>}{|\zeta_n|^2} \\
    &= \left< \frac{\zeta_n}{|\zeta_n|} , -J_0 \cdot [DY_n(\tau_n(t),0) + \epsilon_n] \cdot \frac{\zeta_n}{|\zeta_n|} \right>
  \end{aligned}
\]
where $J_0 = \begin{bmatrix} 0 & -1 \\ 1 & 0 \end{bmatrix}$ and
\[
  \epsilon_n(t) = \int_0^1 [DY_n(\tau_n(t),\lambda\zeta_n(t))-DY_n(\tau_n(t),0)] d\lambda.
\]
Now, since $\hat R_n$ converges $C^1$ to $\hat R_0(\tau,z) = (1,i2z)$, uniformly on $\R/\pi\Z \times \overline{B_\rho(0)}$, we conclude that, as $n\to\infty$, $\epsilon_n$ converges uniformly in $t \in [0,T_n']$ to the zero matrix and $DY_n(\tau_n(t),0)$ converges uniformly in $t \in [0,T_n']$ to the matrix
\[
  \begin{bmatrix} 0 & -2 \\ 2 & 0 \end{bmatrix}
\]
as was computed in~\eqref{DR_n}. Consequently, one estimates $\dot\vartheta_n \geq 3/2$, uniformly in $t\in[0,T_n']$.
Thus $\vartheta_n$ is strictly increasing, and increases at least $3\pi/2$ on each $I^n_j$, for every $j$ and every $n$ sufficiently large. Thus $l_n \to +\infty$, and this contradiction concludes Case 2.

\appendix

\section{Lemmas from Finsler geometry}

\subsection{Proof of Lemma~\ref{useful_lemma}}

For a fixed $j$ we have, in natural coordinates, $(\partial_{x^j})_h = \partial_{x^j} - \Gamma^k_j\partial_{y^k}$. Thus, since $S$ is horizontal, we obtain
\[
  \begin{aligned}
    &\left[ P_H(S), P_H((\partial_{x^j})_h) \right] = [S,(\partial_{x^j})_h] = [y^i\partial_{x^i}-2G^i\partial_{y^i},\partial_{x^j} - \Gamma^k_j\partial_{y^k}] \\
    &= -2[G^i\partial_{y^i},\partial_{x^j}] - [y^i\partial_{x^i},\Gamma^k_j\partial_{y^k}] + 2[G^i\partial_{y^i},\Gamma^k_j\partial_{y^k}] \\
    &= 2(\partial_{x^j}G^i)\partial_{y^i} - y^i(\partial_{x^i}\Gamma^k_j)\partial_{y^k} - \Gamma^k_j[y^i\partial_{x^i},\partial_{y^k}] + 2G^i(\partial_{y^i}\Gamma^k_j)\partial_{y^k} + 2\Gamma^k_j[G^i\partial_{y^i},\partial_{y^k}] \\
    &= 2(\partial_{x^j}G^i)\partial_{y^i} - y^i(\partial_{x^i}\Gamma^k_j)\partial_{y^k} + \Gamma^i_j\partial_{x^i} + 2G^i(\partial_{y^i}\Gamma^k_j)\partial_{y^k} - 2\Gamma^k_j(\partial_{y^k}G^i)\partial_{y^i} \\
    &= \left\{ 2(\partial_{x^j}G^i) - y^k(\partial_{x^k}\Gamma^i_j) + 2G^k(\partial_{y^k}\Gamma^i_j) -2\Gamma^k_j(\partial_{y^k}G^i) \right\} \partial_{y^i} + \Gamma^i_j\partial_{x^i}
  \end{aligned}
\]
and consequently
\begin{equation}\label{gri}
  \begin{aligned}
    & R(S,(\partial_{x^j})_h) = P_V (\left[ P_H(S), P_H((\partial_{x^j})_h) \right]) \\
    &= \left\{ 2(\partial_{x^j}G^i) - y^k(\partial_{x^k}\Gamma^i_j) + 2G^k\Gamma^i_{jk} -2\Gamma^k_j\Gamma^i_k + \Gamma^k_j\Gamma^i_k \right\} \partial_{y^i}
  \end{aligned}
\end{equation}
Here we used $P_V(\partial_{x^i}) = \Gamma_i^k\partial_{y^k}$. On the other hand,
\begin{equation}\label{lift}
  \begin{aligned}
    &\nabla_{[S,(\partial_{x^j})_h]}C = \nabla_{[S,(\partial_{x^j})_h]}(y^l\partial_{y^l}) = (dy^l \cdot [S,(\partial_{x^j})_h])\partial_{y^l} + y^l\nabla_{[S,(\partial_{x^j})_h]}\partial_{y^l} \\
    &= (dy^l \cdot [S,(\partial_{x^j})_h])\partial_{y^l} + y^l\Gamma^i_j\nabla_{\partial_{x^i}}\partial_{y^l} = (dy^l \cdot [S,(\partial_{x^j})_h])\partial_{y^l} + y^l\Gamma^i_j\Gamma^k_{il}\partial_{y^k} \\
    &= \left\{ 2(\partial_{x^j}G^l) - y^k(\partial_{x^k}\Gamma^l_j) + 2G^k\Gamma^l_{jk} -2\Gamma^k_j\Gamma^l_k + y^k\Gamma^i_j\Gamma^l_{ik} \right\} \partial_{y^l}
  \end{aligned}
\end{equation}
Using the identity $y^k\Gamma^l_{ik} = \Gamma^l_i$ one sees that \eqref{gri} equals \eqref{lift}. Since $\nabla_SC=\nabla_{(\partial_{x^j})_h}C=0$ one computes at the base point $v$: $$ \tilde R(v_h,(\partial_{x^j})_h)i_v(v) = \tilde R(S,(\partial_{x^j})_h)C = -\nabla_{[S,(\partial_{x^j})_h]}C = -R(S,(\partial_{x^j})_h). $$ Here we used that the horizontal lift of $v$ to $T_vTM$ is $S$, that $i_v(v) = C$ and, in the last equality, that \eqref{gri} equals \eqref{lift}. The conclusion follows because $i_v$ is an isomorphism.

\subsection{Proof of Lemma~\ref{lemma_metric}}

To prove~\eqref{metric_props} take natural coordinates and write
\[
  \begin{aligned}
    \text{LHS} &= \left( (\partial_{x^l}g_{ij})\dot x^l - 2(\partial_{y^l} g_{ij})G^l - g_{lj}\Gamma^l_i - g_{il}\Gamma^l_j \right) V^iW^j + \text{RHS}.
  \end{aligned}
\]
Since $V$ and $W$ are arbitrary we get
\[
  (\partial_{x^l}g_{ij})\dot x^l - 2(\partial_{y^l} g_{ij})G^l - g_{lj}\Gamma^l_i - g_{il}\Gamma^l_j = 0, \ \forall i,j.
\]
Choosing a symmetric lift $\nabla$ of $\Gamma_S$, with local coefficients $\Gamma^l_{ij}$, the above expression becomes $\left( (\partial_{x^k}g_{ij}) - A_{ijl}(\partial_{y^k}G^l) - g_{lj}\Gamma^l_{ik} - g_{il}\Gamma^l_{jk} \right)\dot x^k = 0$, for every $i,j,k$, where $A_{ijl} = \frac{1}{2}\partial_{y^l} g_{ij}$ are the components of the Cartan tensor~\eqref{cartan_tensor}. Thus, if the symmetric lift satisfies
\begin{equation}\label{chern_pre}
  (\partial_{x^k}g_{ij}) - 2A_{ijl}(\partial_{y^k}G^l) - g_{lj}\Gamma^l_{ik} - g_{il}\Gamma^l_{jk} = 0, \ \forall i,j,k
\end{equation}
the claim follows. However, this may not be satisfied by an arbitrary $\nabla$. In fact, we can permute the $(i,j,k)$ above to get terms corresponding to $(j,k,i)$ and $(k,i,j)$. Adding the $(i,j,k)$-term to the $(j,k,i)$-term, subtracting the $(k,i,j)$-term and using $\Gamma^l_m= \partial_{y^m}G^l$, one obtains
\[
  2g_{jl}\Gamma^l_{ik} = \partial_{x^k}g_{ij} - \partial_{x^j}g_{ki} + \partial_{x^i}g_{jk} - 2\{ A_{ijl}\Gamma^l_k + A_{jkl}\Gamma^l_i  - A_{kil}\Gamma^l_j \}.
\]
Multiplying by $\frac{1}{2}g^{mj}$ and summing in $j$ we get equations~\eqref{chern} for the coefficients of the Chern connection, that is, if we use the Chern connection as the symmetric lift, the desired conclusion holds. However, equation~\eqref{metric_props} does not depend on this choice, which implies that we could have chosen any symmetric lift to carry on our calculations.

\subsection{Proof of Lemma~\ref{lemma_rep_lf}}

Consider, in natural coordinates $x=(x^1,\dots,x^n)$, $y=(y^1,\dots,y^n)$,
the ODE $(\dot x,\dot y) = (y,-2G)$ associated to the geodesic spray, where $G = (G^1,\dots,G^n)$ are the spray coefficients. Linearizing we get
\begin{equation}\label{lin_flow}
  \begin{pmatrix} \dot {\delta x} \\ \dot {\delta y} \end{pmatrix} = \begin{pmatrix} 0 & I \\ -2D_xG & -2D_yG \end{pmatrix} \begin{pmatrix} \delta x \\ \delta y \end{pmatrix}
\end{equation}
where $\delta x = (\delta x^1,\dots,\delta x^n)$ and $\delta y= (\delta y^1,\dots,\delta y^n)$ are the fiber coordinates on $TTM$. Here $D_xG$ and $D_yG$ are $n\times n$ matrices with entries $\partial_{x^i}G^k$ and $\partial_{y^i}G^k$, respectively, evaluated at $(x^1 = x^1\circ\gamma(t),\dots, y^1(t) = \dot x^1,\dots)$, where $\gamma(t)$ is some geodesic. 

We need to rewrite~\eqref{lin_flow} in terms of invariantly defined objects. Consider the vector field $J(t) = \delta x^i \partial_{x^i}$ along $\gamma$. Then $\dot{\delta x^k} = (\frac{D_\gamma J}{dt})^k - \Gamma^k_l\delta x^l$, where $\frac{D_\gamma J}{dt} = (\frac{D_\gamma J}{dt})^k \partial_{x^k}$, and
\begin{equation}\label{delta2}
  \begin{aligned}
    \ddot {\delta x^k} &= \frac{d}{dt}\left[ \left( \frac{D_\gamma J}{dt} \right)^k - \Gamma^k_l\delta x^l \right] = \left( \frac{D^2_\gamma J}{dt^2} \right)^k - \Gamma^k_l  \left( \frac{D_\gamma J}{dt} \right)^l - \frac{d}{dt} (\Gamma^k_l\delta x^l) \\
    &= \left( \frac{D^2_\gamma J}{dt^2} \right)^k - 2\Gamma^k_l \dot{\delta x^l} - (\partial_{x^i}\Gamma^k_l)\dot x^i\delta x^l + 2(\partial_{y^i}\Gamma^k_l)G^i\delta x^l - \Gamma^k_i\Gamma^i_l\delta x^l
  \end{aligned}
\end{equation}
where we used $\dot y^i + 2G^i=0$. Adding $2(\partial_{x^i}G^k)\delta x^i + 2(\partial_{y^i}G^k)\dot{\delta x^i}$ to~\eqref{delta2} we get $0$, in view of~\eqref{lin_flow}, and
\begin{equation}\label{delta3}
  \begin{aligned}
    0 &= \left( \frac{D^2_\gamma J}{dt^2} \right)^k - 2\Gamma^k_l \dot{\delta x^l} - (\partial_{x^i}\Gamma^k_l)\dot x^i\delta x^l + 2(\partial_{y^i}\Gamma^k_l)G^i\delta x^l \\
    &- \Gamma^k_i\Gamma^i_l\delta x^l + 2(\partial_{x^l}G^k)\delta x^l + 2(\partial_{y^l}G^k)\dot{\delta x^l} \\
    &= \left( \frac{D^2_\gamma J}{dt^2} \right)^k + \left\{ 2(\partial_{x^l}G^k) - (\partial_{x^i}\Gamma^k_l)\dot x^i + 2G^i(\partial_{y^i}\Gamma^k_l) - \Gamma^i_l\Gamma^k_i \right\} \delta x^l.
  \end{aligned}
\end{equation}
Now observe that $R^{\dot\gamma}(J) = \delta x^l R^{\dot\gamma}(\partial_{x^l}) = \delta x^l i_{\dot\gamma}^{-1}(-\nabla_{[S,(\partial_{x^l})_h]}C)$ where the vectors $\nabla_{[S,(\partial_{x^l})_h]}C$ were computed in~\eqref{lift}. Plugging into~\eqref{delta3} we get equations
\begin{equation*}\label{}
  \begin{array}{ccc}
    0 = \left( \frac{D^2_\gamma J}{dt^2} \right)^k - R^{\dot\gamma}(J)^k, \ \forall k=1\dots n.
  \end{array}
\end{equation*}

\end{document}